\documentclass[12pt]{paper}
\usepackage{graphicx}
\usepackage{amsmath}
\usepackage{url}
\usepackage{verbatim}
\usepackage{epsf}
\usepackage{latexsym}
\usepackage{array,tabu,multirow,multicol,longtable}
\usepackage{mathrsfs}
\usepackage{authblk}
\usepackage{amsfonts}
\usepackage{amssymb}
\usepackage{moreverb,footmisc,amsthm}
\usepackage{color}
\usepackage[left=3cm,right=2cm,top=2cm,bottom=2cm]{geometry}
\usepackage{hyperref}
\usepackage[toc,page]{appendix}

\hypersetup{
	colorlinks=true,
	linkcolor=blue,
	urlcolor=magenta,
	citecolor=red,
}
\usepackage{parskip}
\usepackage{algorithm}
\usepackage{algorithmicx}
\usepackage{algpseudocode}
\newcommand*{\QED}{\hfill\ensuremath{\square}}%

\newenvironment{Proof}{\vspace{1.5ex}{\sc Proof}. }{\vspace{2ex} $\QED$}

\newtheorem{Theorem}{Theorem}
\newtheorem{Corollary}{Corollary}

\newtheorem{Lemma}{Lemma}
\newtheorem{Proposition}{Proposition}
\newtheorem{Remark}{Remark}

\def\Fc{\mathcal{F}}

\def\Nb{\mathbb{N}}

\def\Rb{\mathbb{R}}

\def\to{\rightarrow}

\def\L{\mathscr{L}}

\def\partder#1#2{
	\frac{\partial #1}{\partial #2}
    }

\def\scal#1#2{
    \left\langle #1,
    #2\right\rangle
    }

\def\mod#1{
    \left|#1
    \right|
    }

\def\norm#1{
    \left\|
    #1\right\|
    }

\def\bigO#1{\mathcal{O}\left(#1\right)}

\def\ud{\,\mathrm{d}}

\begin{document}

\raggedbottom 
\addtolength{\topskip}{0pt plus 10pt}

\title{{\bf \begin{center}
A variational technique of mollification applied to backward heat conduction problems.
\end{center}}}
\author{Walter C. {\sc SIMO TAO LEE}\thanks{Institut de Mathématiques de Toulouse \\
Email: wsimotao@math.univ-toulouse.fr}}
\maketitle

\abstract 

This paper addresses a backward heat conduction problem with fractional Laplacian and time-dependent coefficient in an unbounded domain. The problem models generalized diffusion processes and is well-known to be severely ill-posed. We investigate a simple and powerful variational regularization technique based on mollification. Under classical Sobolev smoothness conditions, we derive order-optimal convergence rates between the exact solution and regularized approximation in the practical case where both the data and the operator are noisy. Moreover, we propose an order-optimal a-posteriori parameter choice rule based on the Morozov principle. Finally, we illustrate the robustness and efficiency of the regularization technique by some numerical examples including image deblurring.
  

\section{Introduction}
This paper deals with the final value problem
\begin{equation}
\label{main equation}
\begin{cases}
\partder{u}{t}(x,t) + \gamma(t) (-\Delta_x)^\tau u(x,t) = 0 & x \in \Rb^n, \,\,\, t \in (0,T)\\
u(x,T) = g(x) & x \in \Rb^n,
\end{cases}
\end{equation}
where we aim at recovering the initial temperature distribution $u(\cdot,0)$ given the time-dependent positive coefficient $\gamma(\cdot)$ and the final temperature distribution $u(\cdot,T)$. In equation \eqref{main equation}, $\tau \in (0,1]$ and $(-\Delta_x)^\tau$ denotes the fractional Laplace operator defined as
$$
\forall \xi \in \Rb^n, \quad \Fc((-\Delta_x)^\tau f)(\xi) = \mod{2\pi\xi}^{2\tau} \Fc(f)(\xi),
$$
where $\Fc(f)$ represents the Fourier transform of a function $f: \Rb^n \to \Rb$ given as
\begin{equation}
\label{def fourier transform}
\forall \xi \in \Rb^n, \quad \Fc(f)(\xi) = \int_{\Rb^n} f(x)e^{- 2\pi i x\cdot\xi} \ud x.
\end{equation}
Of course, for $\tau = 1$, we recover the classical Laplace operator. For $\tau \in (0,1)$, the above definition is equivalent to the definition of fractional Laplacian defined via singular integral operator \cite{kwasnicki2017ten}. 

The dependence in time of the coefficient $\gamma$ indicates that the thermal conductivity varies with time. Such is the case of heat conduction in a material subject to radio-active decay, the thermal conductivity being dependent of the amount of decay is then related to time. The fields of application of backward heat conduction problems are quite diversified. Beyond applications in heat transfer and diffusion problems \cite{engl1995inverse} where $\tau$ is usually equal to $1$, we can also speak of deconvolution problems and deblurring processes \cite{carasso1994overcoming,carasso2003apex,carasso1978digital} where $\tau \in (0,1)$ and hydrology \cite{skaggs1994recovering}. 

   Backward heat conduction problems are well-known to be ill-posed \cite{isakov2006inverse,payne1975improperly} in the sense that small perturbation of the final distribution of temperature $u(\cdot, T)$ may induce arbitrary large errors in the initial distribution of temperature $u(\cdot, 0)$.
Given the ill-posedness of the problem, the use of a regularization method is critical in order to recover stable approximate of initial distribution of temperature. In this regard, several regularization techniques have been proposed in the literature for solving final value heat conduction problems. We can mention, the method of fundamental solutions \cite{mera2005method}, methods using finite difference approximation \cite{elden1995numerical,iijima2004lattice,shidfar2005numerical}, the method of quasi-reversibility \cite{methodequasireversibiliteLions,miller1973stabilized,showalter1974final}, iterative methods \cite{jourhmane2002iterative,kozlov1989iterative,mera2001iterative}, truncation methods \cite{tuan2012determination,tuan2011two,zhang2014posteriori}, optimal filtering \cite{seidman1996optimal}, operator-splitting method \cite{kirkup2002solution}, kernel-based method \cite{ames1997kernel}, wavelet regularization \cite{wang2011shannon}, modified method \cite{le2013backward,qian2007modified}, Fourier regularization \cite{fu2007fourier,fu2012posteriori}, quasi boundary regularization \cite{chang2007new,denche2005modified}, non local boundary value problem method \cite{hao2011stability,hao2008non} and Tikhonov method \cite{hon2011discretized,ji2002numerical}. Besides those methods, there are also methods based on mollification \cite{h1994mollification,hao2009stability,van2017posteriori}.

Despite the considerable amount of papers treating backward heat conduction problems, few treat the case on general dimension $n \in \Nb$ with time-dependent conductivity coefficient $\gamma$. Moreover, among those who did so, it is very seldom to see a paper that covers the case of fractional Laplacian. 
This article presents a variational technique based on mollification for the regularization of backward heat equation with time-dependent coefficient and fractional Laplician in the unbounded domain $\Rb^n$, $n \in \Nb$. We recall that there are three main mollification approaches for regularizing an ill-posed equation 
\begin{equation*}
A u = g,
\end{equation*}
where $A: U \to G$ is a bounded operator between Banach spaces $U$ and $G$, $g$ is the data and $u$ is the unknown solution.
     One technique consists in pre-smoothing the data $g$ by a mollifier operator $M_\alpha$ and check that the equation $ A u = M_\alpha g$ is well-posed, in which case, the sought solution $u$ is approximated by $u_\alpha = A^\dagger M_\alpha g$ which is nothing but the solution of the regularized equation $ A u = M_\alpha g$. This approach has been studied by Vasin \cite{vasin1973stable}, Murio \cite{murio1987automatic,murio2011mollification}, Manselli and Miller \cite{manselli1980calculation}, H\`ao et al. \cite{hao1994stable,hao2009stability,h1994mollification}, Van Duc \cite{van2017posteriori}. Notice that this approach usually requires the application of Moore-Penrose pseudo inverse $A^\dagger$ of $A$ which is unbounded, though the pre-composition by the mollifier operator $M_\alpha$ (i.e. $A^\dagger M_\alpha$) is bounded. In \cite{hao2009stability,h1994mollification,van2017posteriori}, this approach has been applied to backward heat equation.
     A second approach based on the adjoint equation is the so called \textit{approximate-inverse} developped by Louis and Mass \cite{louis1996approximate,louis1999unified,louis1990mollifier}. In this case, $U$ and $G$ are functional Hilbert spaces, and the solution $u$ is approximated by a mollifier version $E_\gamma u$ defined using inner product as $ (E_\gamma u)(x) = \left\langle e_\gamma(x,\cdot) ,u\right\rangle$ where $e_\gamma(x,\cdot)$ lives in the range of the adjoint operator $A^*$ of $A$, i.e. $e_\gamma(x,\cdot) = A^* v_{x,\gamma}$ for some $v_{x,\gamma} \in G$. Then, using the adjoint equation, the approximate solution $(E_\gamma u)$ is computed as $(E_\gamma u)(x) = \left\langle  v_{x,\gamma} , g  \right\rangle.$ Despite many advantages of this approach, the main limitation is the fact that such explicit pair of \textit{reconstruction kernel - mollifier} $(v_{x,\gamma},e_\gamma)$ is known only in few cases. Otherwise, if such pair is unknown, given $e_\gamma$, one has to solves the adjoint equation $A^* v_{x,\gamma} = e_\gamma(x,\cdot)$ (which is ill-posed though with exact data) in order to deduce the \textit{reconstruction kernel} $v_{x,\gamma}$ which allows to compute the approximate solution $E_\gamma u$. A good account of this method and some applications is given in the book by Schuster \cite{schuster2007method}.
	The third approach of mollification is a variational one, which to the best knowledge of the author is due to Lannes et al. \cite{lannes1987stabilized} in signal and image processing. This way of regularization was then further studied by Alibaud et al. \cite{alibaud2009variational} where the authors proved for the first time consistency results of the method, in the particular problem of Fourier synthesis. In the main-time, Bonnefond and Maréchal \cite{bonnefond2009variational} proposed a generalization of the method to the inversion of some compact operators. However, despite the flexibility and interesting applications of this approach (see e.g. \cite{marchal2000new}), no error estimates or rates of convergence of the method appears in the literature. The best results available deal only with consistency of the method. 
	
In this paper, we propose and study a variational approach of mollification closely related to the one studied by Alibaud et al. \cite{alibaud2009variational} and Bonnefond and Maréchal \cite{bonnefond2009variational}. The aim of this paper is to present a simple and flexible variational regularization by mollification, with not only consistency results but also order-otimal convergence rates as well as a-posteriori parameter selection rules leading to order-optimal rates. More precisely, we first formulate equation \eqref{main equation} in the form of an operator equation 
\begin{equation}
\label{op eq of our pb}
A u(\cdot,0) = g,
\end{equation}
with $A: L^2(\Rb^n) \to L^2(\Rb^n)$ is bounded. Then we define the regularized solution. Under noisy data $g^\delta$ and approximate operator $A^h$ verifying 
\begin{equation*}
||g^\delta - g||_{L^2} \leq \delta \quad \text{and} \quad \vert|| A - A^h \vert|| \leq h,
\end{equation*}
we derive order-optimal convergence rates under the classical Sobolev smoothness condition
\begin{equation}
\label{smoothness cond}
u(\cdot,0) \in H^p(\Rb^n) \quad \text{with} \quad \norm{u(\cdot,0)}_{H^p} \leq E. 
\end{equation}
Moreover, we show that the smoothness condition \eqref{smoothness cond} is equivalent to some {\it logarithmic source condition}  which is the natural type of smoothness condition in context of exponentially ill-posed problems \cite{hohage2000regularization,tautenhahn1998optimality}. 

The outline of this article is as follows: 
In Section \ref{section regularization}, We reformulate the partial differential equation \eqref{main equation} into an operator equation of the form \eqref{op eq of our pb} and define the variational regularization technique together with some interesting remarks regarding the penalty term in comparison to classical variational approaches. Section \ref{section error estimate} deals with error estimates and order-optimality of the method under the smoothness condition \eqref{smoothness cond} for both the cases of exact or approximately given operator $A$. 
Section \ref{section par choice rule} is devoted to parameter selection rules which is a critical step in the application of a regularization method. Here we propose a Morozov-like a-posteriori parameter choice rule leading to order-optimal convergence rates. Finally, we study four numerical examples, including an image deblurring, in Section \ref{section numerical experiments} which illustrate the great potential and attractiveness of the regularization approach. Moreover, in this Section, the numerical convergence rates observed confirm the theoretical convergence rates given in the paper.

We point out that all the results in this paper cover the classical backward heat equation with Laplace operator which is just the particular case $\tau = 1$. In the sequel, $\hat{f}$ or $\Fc(f)$ (resp. $\Fc^{-1}(f)$ ) denotes the Fourier transform (resp. the inverse Fourier transform) of the function $f$ defined by \eqref{def fourier transform}, $||f||$ or $||f||_{L^2}$ always refers to the $L^2$-norm of the function $f$ on $\Rb^n$, $||f||_{H^p}$ denotes the Sobolev norm of $f$ on $\Rb^n$ and $||| \cdot ||||$ denotes operator norm of a bounded linear mapping. 
\section{Regularization}\label{section regularization}

Before goint to the regularization procedure of the partial differential equation \eqref{main equation}, let us first put it into a well defined framework. We recall that we aim to recover the initial temperature $u(\cdot,0)$ from the the final temperature distribution $g = u(\cdot,T)$. Henceforth, we assume that the data $g$ belongs to $L^2(\Rb^n)$.

By applying the Fourier transform in \eqref{main equation}, we readily get
\begin{equation}
\label{fourier equation}
\begin{cases}
\partder{\hat{u}}{t}(\xi,t) = - \gamma(t) |2\pi\xi|^{2\tau} \hat{u}(\xi,t) & \xi \in \Rb^n, \,\,\, t \in (0,T)\\
\hat{u}(\xi,T) = \hat{g}(\xi) & \xi \in \Rb^n,
\end{cases}
\end{equation}
from which we derive that
\begin{equation}
\label{link equation data and sol in fourier domain}
\forall \xi \in \Rb^n,\quad \hat{u}(\xi,T) = \hat{u}(\xi,0) \exp\left(-  | 2\pi\xi|^{2\tau} \int_0^T \gamma(\lambda)\ud \lambda\right).
\end{equation}
In \eqref{link equation data and sol in fourier domain}, we assume that the conductivity function $\gamma$ is positive and integrable on $(0,T)$.
Hence, using the Fourier transform, we can rewrite problem \eqref{main equation} into an operator equation 
\begin{equation}
\label{operator equation problem}
A u(\cdot,0) = g
\end{equation}
where $A: L^2(\Rb^n) \to L^2(\Rb^n)$ is defined by
\begin{equation}
\label{def operator A}
\forall f \in L^2(\Rb^n), \quad A f = \Fc^{-1} \left( \exp\left(-|2\pi\xi|^{2\tau} \int_0^T \gamma(\lambda)\ud \lambda\right) \hat{f}(\xi) \right),
\end{equation}
from which we formally write
$$
A = \Fc^{-1}  \psi(\xi)  \Fc,
$$
with the function $\psi$ defined by
\begin{equation}
\label{def funct psi}
\psi(\xi) = \exp\left(-|2\pi\xi|^{2\tau} \int_0^T \gamma(\lambda)\ud \lambda\right).
\end{equation}
Using the Parseval identity and the positivity of the function $\gamma$, we straightforward get the boundedness of operator $A$ with $|||A||| \leq 1$. Also the exponential ill-posedness of equation \eqref{operator equation problem} stems from
$$
\hat{u}(\cdot,0) = \hat{g}(\xi) \exp\left(|2\pi \xi|^{2\tau} \int_0^T \gamma(\lambda)\ud \lambda\right),
$$
which leads to unbounded exponential amplification of large frequency components in the data $g$. We point out that from \eqref{def operator A}, we can easily notice that $A$ is actually a convolution operator by the inverse Fourier transform $\Fc^{-1}( \psi(\xi))$ of the function $\psi$.

From \eqref{operator equation problem} and \eqref{def operator A}, we set problem \eqref{main equation} into a well-defined framework where given an $L^2$ data $g$, we seek $L^2$ solution $u(\cdot,0)$ satisfying equation \eqref{fourier equation}.

\textbf{Regularized solution}

Given the ill-posedness of equation \eqref{operator equation problem}, it is hopeless to recover $u(\cdot,0)$ from $g$ without the use of a regularization method. Let $\phi$ be a smooth real-valued function in $L^1(\Rb^n)$ satisfying $\int_{\Rb^n} \phi(x) \ud x =1$. It is well-known that the family of functions $(\phi_\beta)_{\beta>0}$ defined by 
\begin{equation}
\label{def phi beta from phi}
\forall x\in \Rb^n,\quad \phi_\beta(x) := \frac{1}{\beta^n}\phi(x/\beta),
\end{equation}
 satisfies
\begin{equation}
\label{prop phi beta}
\forall f \in L^2(\Rb^n), \quad \phi_\beta \star f \to f \quad \text{as} \quad \beta \downarrow 0,
\end{equation}
where $f \star  g$ is the convolution of the functions $f$ and $g$ defined as
$
\left( f\star g \right)(x) = \int_{\Rb^n} f(x-y)g(y) \ud y.
$
For $\beta>0$, let $C_\beta$ be the mollifier operator defined by
\begin{equation}
\label{def oper Cbeta}
\forall \beta > 0, \,\, \forall f \in L^2(\Rb^n), \quad C_\beta f = \phi_\beta \star f.
\end{equation}
From \eqref{prop phi beta}, we see that the family of operators $(C_\beta)_{\beta>0}$ is an approximation of unity  in $L^2(\Rb^n)$, that is,
\begin{equation*}
\label{conver Cbeta}
\forall f \in L^2(\Rb^n), \quad C_\beta f \to f \quad \text{in} \quad L^2(\Rb^n) \quad \text{as} \quad \beta\downarrow 0 .
\end{equation*}

Now, given a data $g \in L^2(\Rb^n)$, we define the regularized solution $u_\beta$ as
\begin{equation}
\label{def reg solution}
u_\beta = \underset{u \in L^2(\Rb^n)}{\mathrm{argmin}} \quad || A u - g||_{L^2}^2 + || (I - C_\beta)u||_{L^2}^2
\end{equation}
where $I$ denotes the identity operator on $L^2(\Rb^n)$.
For convenience, we call $J_\beta(u,A,g)$ the functional minimised in \eqref{def reg solution}, i.e.
\begin{equation}
\label{def func J beta}
J_\beta(u,A,g) := || A u - g||_{L^2}^2 + || (I - C_\beta)u||_{L^2}^2.
\end{equation}
Before proving that $u_\beta$ is well defined and that it induces a regularization method, let us comment a bit on \eqref{def reg solution}.
In the functional $J_\beta(\cdot,A,g)$ minimized in \eqref{def reg solution}, we have a fit term $|| A u - g||^2$ which aims at fitting the original model \eqref{operator equation problem}, a penalty term $|| (I - C_\beta)u||^2$ which aims to introduce stability in the model and a regularization parameter $\beta>0$ which allows to control the level of regularization. Regarding the penalty term, we see that we actually minimize the $L^2$ distance between $u$ and its mollifier version $C_\beta u$. Heuristically, we aim at recovering a smooth version of the unknown solution $u(\cdot,0)$ which is more tractable and yet not far from the true solution itself. Notice that on the contrary to the variational approach of Alibaud et al. \cite{alibaud2009variational} and Bonnefond and Maréchal \cite{bonnefond2009variational}, we do not put a fitting operator on the data $g$.

Now let us state the following theorem which asserts the well-posedness of $u_\beta$ in \eqref{def reg solution} and the fact that this actually induces a regularization method.
\begin{Theorem}
\label{Theorem existence, uniqueness and regularization}
For all data $g \in L^2(\Rb^n)$, the regularized solution $u_\beta$ expressed by \eqref{def reg solution} is well-defined and is characterized by the equation
\begin{equation}
\label{char u beta}
[ A^*A + (I-C_\beta)^*(I-C_\beta) ]\, u_\beta = A^* g .
\end{equation}
 Moreover the mapping $R_\beta: L^2(\Rb^n) \to L^2(\Rb^n)$, $g \mapsto u_\beta$ is linear and bounded and 
\begin{equation}
\label{consistency result}
\forall u \in L^2(\Rb^n), \quad R_\beta A u \to u \quad \text{in} \quad L^2(\Rb^n) \quad \text{as} \quad \beta \downarrow 0.
\end{equation}
That is, \eqref{def reg solution} defines a regularization method for equation \eqref{operator equation problem}.
\end{Theorem} 

\begin{Proof}
For $u \in L^2(\Rb^n)$, we have
\begin{eqnarray}
\label{eq1}
|| A u||^2 + || (I - C_\beta)u||^2 & = & \scal{A^*A + (I - C_\beta)^*(I - C_\beta) u}{u} \nonumber \\
& =& \int_{\Rb^n} \left( \psi(\xi)^2 + |1-\Fc(\phi_\beta)(\xi)|^2 \right) |\hat{u}(\xi)|^2\ud \xi \nonumber \\
& =& \int_{\Rb^n} \left( \psi(\xi)^2 + |1-\hat{\phi}(\beta \xi)|^2 \right) |\hat{u}(\xi)|^2\ud \xi  \nonumber  \\
& \geq & \mu_\beta ||u||^2,
\end{eqnarray} 
where 
$$
\mu_\beta = \min_{\xi \in \Rb^n} \quad \Theta_\beta(\xi) := \psi(\xi)^2 + |1-\hat{\phi}(\beta \xi)|^2.
$$
The existence of $\mu_\beta >0$ is justified by the fact that the function  $\Theta_\beta$ is continuous, strictly positive and tends to $1$ as $\xi \to \infty$.
Hence from \eqref{eq1}, and the inequality $||Au -g||^2 + ||g||^2 \geq ||Au||^2/2$ we deduce that the mapping $u \mapsto J_\beta(u,A,g)$ is coercive as
$$
J_\beta(u,A,g) \geq \frac{\mu_\beta}{2} ||u||^2 - ||g||^2.
$$
In summary, the mapping $u \mapsto J_\beta(u,A,g)$ is well defined, lower semi-continuous, coercive and strictly convex whence the existence and uniqueness of its minimizer $u_\beta$. The characterization \eqref{char u beta} follows readily by application of the first order optimality condition in \eqref{def reg solution}.
Notice also that \eqref{eq1} implies that the operator $ A^*A + (I - C_\beta)^*(I - C_\beta)$ is continuously invertible and that 
$$
|||\left[A^*A + (I - C_\beta)^*(I - C_\beta) \right]^{-1} ||| \leq 1/\mu_\beta.
$$
Taking the Fourier transform in \eqref{char u beta}, we get
\begin{equation}
\label{eq2}
\left( \psi(\xi)^2 + |1-\hat{\phi}(\beta \xi)|^2 \right) \widehat{u_\beta}(\xi) = \psi(\xi)\hat{g}(\xi),
\end{equation}
which implies
\begin{equation}
\label{expr u beta hat}
|\widehat{u_\beta}(\xi) | = \frac{\psi(\xi)}{\psi(\xi)^2 + |1-\hat{\phi}(\beta \xi)|^2} |\hat{g}(\xi) | \leq \frac{1}{\mu_\beta}|\hat{g}(\xi) |.
\end{equation}
The Parseval identity and \eqref{expr u beta hat} implies the boundedness of the mapping $R_\beta$ with~$|||R_\beta||| \leq 1/\mu_\beta$. \\
Let $g = A u$ with $u \in L^2(\Rb^n)$ and $u_\beta = R_\beta g$, then from \eqref{eq2}, we get that
$$
\left( \psi(\xi)^2 + |1-\hat{\phi}(\beta \xi)|^2 \right) \widehat{u_\beta}(\xi) = \psi(\xi)^2\hat{u}(\xi),
$$
which yields
\begin{equation}
\label{norm residual}
|| u_\beta - u||^2 = || \widehat{u_\beta} - \hat{u}||^2  = \int_{\Rb^n} \mod{\frac{|1-\hat{\phi}(\beta \xi)|^2}{\psi(\xi)^2 + |1-\hat{\phi}(\beta \xi)|^2}}^2 |\hat{u}(\xi)|^2 \ud \xi \leq \int_{\Rb^n}  |\hat{u}(\xi)|^2 \ud \xi.
\end{equation}
Given that $\hat{\phi}(\beta \xi) \underset{\beta \to 0}{\to} \hat{\phi}(0) = \int_{\Rb^n} \phi(x) \ud x = 1$, by applying dominated convergence theorem to \eqref{norm residual}, we deduce that $u_\beta$ converges to $u$ in $L^2(\Rb^n)$ as $\beta \downarrow 0$.
\end{Proof}

Now in order to derive error estimates between exact solution $u(\cdot,0)$ and regularized solution $u_\beta$, we need to assess the speed of convergence of the penalty term $||(I -C_\beta)u||^2$ as $\beta$ goes to $0$. Notice that this is generally trivial for classical variational regularization methods such as Tikhonov method, in which case, the regularization parameter, generally named $\alpha$, is a weight in the penalty term. In this regard, let us state the following proposition (whose proof is deferred to appendix) which says that without imposing smoothness condition on the kernel function $\phi$, operator $C_\beta$ converges arbitrarily slowly to the identity operator in $L^2(\Rb^n)$.
\begin{Proposition}
\label{Prop arbitrary slow convergence of C_beta}
Let the family $(C_\beta)_{\beta>0}$ be defined by \eqref{def oper Cbeta} with $\phi \in L^1(\Rb^n)$ satisfying $\int_{\Rb^n} \phi(x) \ud x = 1$. Then the point-wise convergence of $C_\beta$ toward the identity operator $I$ on $L^2(\Rb^n)$ is arbitrarily slow. That is, there does not exist a non-decreasing function  $\nu: \Rb^+ \to \Rb^+$ such that $\lim_{\beta \downarrow 0} \nu(\beta) = 0$ and
\begin{equation}
\label{unif bound I-Cbeta}
 \forall \, f \in L^2(\Rb^n), \quad \norm{(I-C_\beta)f}_{L^2} \leq \nu(\beta)\, \norm{f}_{L^2}.
\end{equation} 
\end{Proposition}
Given Proposition \ref{Prop arbitrary slow convergence of C_beta}, we need to impose an additional assumption on the function $\phi$. 
In the sequel, we assume that the function $\widehat{\phi}$ is radially decreasing and satisfies 
\begin{equation}
\label{cond on conv kernel phi}
|1-\widehat{\phi}(\xi)| \sim \vert\xi\vert^s, \,\,\text{as} \,\, \xi \to 0, \quad \text{with} \quad s>0.
\end{equation}
With condition \eqref{cond on conv kernel phi}, we can establish speed of convergence of the penalty term $\norm{(I - C_\beta)u}$ when $u$ belongs to Sobolev spaces $H^p(\Rb^n)$.

\begin{Lemma}
\label{lemma Alibaud et al}
Let $\phi$ be a function in $L^1(\Rb^n)$ satisfying \eqref{cond on conv kernel phi}  and such that $\int_{\Rb^n} \phi(x) \ud x = 1$ with $|\hat{\phi}(\xi)| < 1$ for $\xi \neq 0$. Let
\begin{equation}
\label{def m_beta and M_beta}
m_\beta = \min_{\vert\xi\vert=1} |1 - \hat{\phi}(\beta \xi)|^2,\quad  \text{and} \quad  M_\beta = \max_{\vert\xi\vert=1} |1 - \hat{\phi}(\beta \xi)|^2.
\end{equation}
Then the following hold:
\begin{itemize}
\item[i)] $ 0< m_\beta \leq M_\beta \leq (1+||\phi||_{L^1})^2$,

\item[ii)] $M_\beta \to 0 \,\, \text{as}\,\, \beta \to 0$ and $\sup_{\beta \in (0,1]} \frac{M_\beta}{m_\beta} < \infty$,

\item[iii)] there exist positive constants $\nu_0$ and $C_0$ such that, for all $\beta \in (0,1]$ and every $\xi \in \Rb^n \setminus \{0\}$,
\begin{equation}
\label{key estimate alibaud}
\nu_0 \left(\vert \xi \vert^{2s} 1_{\{\vert \xi \vert \leq 1/\beta\}} + \frac{1}{M_\beta} 1_{\{\vert \xi \vert > 1/\beta\}} \right) \leq  \frac{|1-\hat{\phi}(\beta \xi)|^2}{|1-\hat{\phi}(\beta \xi/\vert \xi \vert)|^2} \leq C_0 \vert \xi \vert^{2s}.
\end{equation}
\end{itemize}
\end{Lemma}
This lemma can be found together with its proof in \cite{alibaud2009variational}. 
The next lemma completes Lemma \ref{lemma Alibaud et al} and gives speed of convergence of operator $I-C_\beta$ when restricted to sobolev subspaces $H^p(\Rb^n)$.
\begin{Lemma}
\label{lemma completion lemma alibaud}
Consider the setting of Lemma \ref{lemma Alibaud et al} and let $p>0$. Then the following hold:
\begin{itemize}
\item[i)]
$
m_\beta \sim \beta^{2s} \quad \textrm{and} \quad M_\beta \sim \beta^{2s} \quad \textrm{as} \quad \beta \downarrow 0.
$
\item [ii)] There exists a constant $C_1>0$ dependent on $p$ such that
\begin{equation}
\label{upp bound norm f_beta}
\forall f \in H^p(\Rb^n), \quad || (I - C_\beta)f||_{L^2}^2 \leq C_1 \beta^{2(p \wedge s)} ||f||_{H^p}^2,
\end{equation}
\item[iii)] There exists a constant $C_2>0$ dependent on $p$ such that
\begin{equation}
\label{upp bound (I-C_beta)^2 f}
\forall f \in H^{2p}(\Rb^n), \quad || (I - C_\beta)^*(I - C_\beta)f||_{L^2}^2   \leq C_2 \beta^{4(p \wedge s)} ||f||_{H^{2p}(\Rb^n)}^2,
\end{equation}
\end{itemize}
In \eqref{upp bound norm f_beta} and \eqref{upp bound (I-C_beta)^2 f}, $p\wedge s = \min \{p,s\}$.
\end{Lemma}
The proof of this Lemma is deferred to appendix.
With Lemmas \ref{lemma Alibaud et al} and \ref{lemma completion lemma alibaud}, we are ready for error estimates analysis.

\section{Error estimates}\label{section error estimate}

Henceforth, $g^\delta \in L^2(\Rb^n)$ denotes a noisy data satysfing the noise level condition
\begin{equation}
\label{noise level cond on data}
|| g - g^\delta || \leq \delta,
\end{equation}
where $g = u(\cdot,T)$ is the exact final distribution of temperature. The regularized solution $u_\beta^\delta$ corresponding to the noisy data $g^\delta$ is defined as
\begin{equation}
\label{def reg sol noisy data}
u_\beta^\delta = \mathrm{argmin}_{u \in L^2(\Rb^n)} \quad J_\beta(u,A,g^\delta),
\end{equation} 
with $J_\beta$ defined in \eqref{def func J beta}.
Let us state the next proposition which provides estimates of the propagated data noise error $u_\beta - u_\beta^\delta$ between the regularized solution corresponding to exact and noisy data.

\begin{Proposition}
\label{Prop estimate data propagated error}
Consider the setting of Lemma \ref{lemma Alibaud et al}. Let $g^\delta$ be a noisy data satisfying \eqref{noise level cond on data}with $g = u(\cdot,T)$ and $u_\beta$ resp. $u_\beta^\delta$ be defined by \eqref{def reg solution} (resp. \eqref{def reg sol noisy data}), then there exists a constant $C >0$ independent of $\delta$ and $\beta$ such that
\begin{equation}
\label{bound error prog term}
||u_\beta - u_\beta^\delta||_{L^2} \leq C \frac{\delta}{\beta^{s}} \quad \text{as} \quad \beta \downarrow 0.
\end{equation}
\end{Proposition}

\begin{Proof}
Let $\beta>0$. By linearity of the mapping $R_\beta: L^2(\Rb^n) \to L^2(\Rb^n), g \mapsto \mathrm{argmin}_{u \in L^2(\Rb^n)} \,\, J_\beta(u,A,g)$, we get 
$$
u_\beta - u_\beta^\delta = \mathrm{argmin}_{u \in L^2(\Rb^n)} J_\beta(u,A,g-g^\delta),
$$
which yields using Parseval identity
\begin{eqnarray}
\label{eq3}
||u_\beta - u_\beta^\delta||_{L^2(\Rb^n} & = & \int_{\Rb^n)}\mod{ \frac{|\psi(\xi)|}{|\psi(\xi)|^2 + |1- \hat{\phi}(\beta \xi)|^2}}^2 |\hat{g}(\xi) - \widehat{g^\delta}(\xi)|^2 \ud \xi \nonumber \\
 & =  & \int_{|\xi|\leq r} \Psi_\beta(\xi)^2\,|\hat{g}(\xi) - \widehat{g^\delta}(\xi) |^2 \ud\xi + \int_{|\xi| > r} \Psi_\beta(\xi)^2\,|\hat{g}(\xi) - \widehat{g^\delta}(\xi) |^2 \ud\xi ,
\end{eqnarray}
where $r$ is a positive number and $\Psi_\beta$ is defined by
\begin{equation*}
\label{def func Psi_beta moll deconv}
\Psi_\beta(\xi) = \frac{|\psi(\xi)|}{|\psi(\xi)|^2 + |1-\hat{\phi}(\beta\xi)|^2}.
\end{equation*}
Given that $\psi$ is radially decreasing function, there exists $r_0>0$ such that $|\psi(\xi)| > \frac{1}{2}|\psi(0)|$ for all $|\xi| \leq r_0$.
For $|\xi| \leq r_0$, we have $\Psi_\beta(\xi) \leq 1/|\psi(\xi)| \leq 2/|\psi(0)|$ which implies that
\begin{equation}
\label{inte mol conv}
\int_{|\xi|\leq r_0} \Psi_\beta(\xi)^2\,|\hat{g}(\xi) - \widehat{g^\delta}(\xi) |^2 \ud\xi  \leq  \frac{4}{|\hat{\psi}(0)|^2} \int_{|\xi|\leq r_0} |\hat{g}(\xi) - \widehat{g^\delta}(\xi) |^2 \ud\xi .
\end{equation}
Let $h_0(\beta) = \min \left\lbrace |1 - \hat{\phi}(\beta \xi )|, \quad \xi \in \Rb^n \,\,\, \text{s.t.} \,\,\, |\xi|\geq r_0 \right\rbrace$. Given the continuity of the function $\hat{\phi}$, the fact that $|1 - \hat{\phi}(\beta \xi )| \to 1$ as $|\xi| \to \infty$ and the fact that $\widehat{\phi}$ is radially decreasing, the minimum $h_0(\beta)$ exists and is attained at some $\xi_0 \in \Rb^n$ satisfying $|\xi_0| = r_0$. Moreover, $h_0(\beta)>0$ as $|\hat{\phi}(\xi)| < 1$ for $\xi \neq 0$. So, we deduce that
\begin{eqnarray*}
\forall \,|\xi|\geq r_0, \quad\Psi_\beta(\xi) \leq \frac{|\hat{\psi(\xi)}|}{2 |\hat{\psi}(\xi)||1-\hat{\phi}(\beta\xi)|} \leq \frac{1}{2 |1-\hat{\phi}(\beta \xi_0)|}.
\end{eqnarray*}
This implies that 
\begin{equation}
\label{inte 2 mol deconv}
\int_{|\xi| > r_0} \Psi_\beta(\xi)^2\,|\hat{g}(\xi) - \widehat{g^\delta}(\xi) |^2 \ud\xi  \leq  \frac{1}{4 |1-\hat{\phi}(\beta \xi_0)|^2} \int_{|\xi| > r} |\hat{g}(\xi) - \widehat{g^\delta}(\xi) |^2 \ud\xi.
\end{equation}
Applying \eqref{eq3} with $r =r_0$ together with \eqref{inte mol conv} and \eqref{inte 2 mol deconv} yields
\begin{equation}
\label{yyyy mol deconv}
||u^\beta - u_\beta^\delta||^2 \leq \max\left\lbrace\frac{4}{|\hat{\psi}(0)|^2}, \frac{1}{4 |1-\hat{\phi}(\beta \xi_0)|^2} \right\rbrace ||\hat{g}- \widehat{g^\delta} ||^2.
\end{equation}
Since $|1-\hat{\phi}(\beta \xi_0)| \to 0$ as $\beta$ goes to $0$, we deduce that for small $\beta$, $ |\hat{\psi}(0)|^2/4 > 4 |1-\hat{\phi}(\beta \xi_0)|^2$ . Then \eqref{yyyy mol deconv} implies that
\begin{equation}
\label{zzz mol deconv}
||u_\beta - u_\beta^\delta|| \leq  \frac{||\hat{g}- \widehat{g^\delta} ||.}{ 2|1-\hat{\phi}(\beta \xi_0)|} \quad \text{as} \quad \beta \downarrow 0. 
\end{equation}
The estimate \eqref{bound error prog term} follows immediately by applying \eqref{cond on conv kernel phi} and \eqref{noise level cond on data} to \eqref{zzz mol deconv}.
\end{Proof}

From \eqref{consistency result} and \eqref{bound error prog term}, we can deduce that for $\beta(\delta) = \delta^{\theta/s}$ with $\theta \in (0,1)$, the reconstruction error between the exact solution $u(\cdot,0)$ and the final reconstruction $u_\beta^\delta$ converges to $0$ as the noise level $\delta$ goes to $0$. That is, if $\beta(\delta) = \delta^{\theta/s}$ with $\theta \in (0,1)$, then
\begin{equation}
\label{conver method noisy case}
||u(\cdot,0) - u_{\beta(\delta)}^\delta|| \to 0 \quad \text{as} \quad \delta \to 0.
\end{equation}
 
Now let us study the regularization error $u(\cdot,0) - u_\beta$. It is well known   that without imposing smoothness condition on the exact solution $u(\cdot,0)$ (or on the exact data $g$), the regularization error of any regularization method converges arbitrarily slowly to $0$ (see, e.g. \cite{schock1985approximate}). Henceforth, we consider the smoothness condition 
\begin{equation}
\label{smoothness cond on u(cdot,0)}
\norm{u(\cdot,0)}_{H^p} \leq E, \quad \text{with} \quad  p>0, \,\,E>0.
\end{equation}
Before getting into the analysis of convergence rates, let us first see how the smoothness \eqref{smoothness cond on u(cdot,0)} is linked to logarithmic source conditions generally occurring in the regularization of exponentially ill-posed problems.

Henceforth, $\bar{A} := A/\sqrt{e}$ denotes the normalized operator $A$, that is, $ \bar{A} = \Fc^{-1}  \frac{\psi(\xi)}{\sqrt{e}}  \Fc$. In fact, with this normalization, we have $|||\bar{A}^*\bar{A}||| \leq e^{-1}$ which is necessary for logarithmic source condition (see, e.g. \cite{hohage2000regularization}).
\begin{Proposition}
\label{Prop 1 smoothness condition}
There exists a constant $C(p)$ depending on $p$ such that the smoothness condition \eqref{smoothness cond on u(cdot,0)} is equivalent to the logarithmic source condition 
\begin{equation}
\label{log sour cond}
u(\cdot,0) = f_{q}(\bar{A}^*\bar{A}) w \quad \text{with} \quad ||w||_{L^2} \leq \rho \quad \text{with} \quad \rho = C(p) E,
\end{equation}
where $q = p/2\tau$ and $f_{q}$ is defined as
\begin{equation*}
\forall \lambda \in (0,1), \quad f_{q}(\lambda) = (-\ln(\lambda))^{-q}.
\end{equation*}
\end{Proposition}

\begin{Proof}
$\Rightarrow$) Assume \eqref{smoothness cond on u(cdot,0)}. Consider the function $\varpi$ defined by 
\begin{equation}
\label{def varpi}
\varpi(\xi)  = \left( 1 + 2|2\pi\xi|^{2 \tau} \int_0^T \gamma(\lambda) \ud \lambda \right)^{p/2 \tau} \hat{u}(\xi,0).
\end{equation}
We have
\begin{equation}
\label{eq4}
|\varpi(\xi)| \leq \max \left\lbrace 1, \left( 2 (2\pi)^{2 \tau} \int_0^T  \gamma(\lambda) \ud \lambda \right)^{p/2\tau} \right\rbrace \left(1 + |\xi|^{2\tau} \right)^{p/2\tau} |\hat{u}(\xi,0)|.
\end{equation}
But since $u(\cdot,0) \in H^p(\Rb^n)$ then $\left(1 + |\xi|^{2\tau} \right)^{p/2\tau} |\hat{u}(\xi,0)| \in L^2(\Rb^n)$ so that \eqref{eq4} implies that $\varpi \in L^2(\Rb^n)$.
By letting $ w =\Fc^{-1} \varpi$, \eqref{def varpi} implies that
\begin{equation}
\label{eq5}
\hat{u}(\xi,0) = \left( 1 + 2 |2\pi\xi|^{2 \tau} \int_0^T  \gamma(\lambda) \ud \lambda \right)^{-p/2 \tau} \hat{w}(\xi)
\end{equation}
By noticing that $\bar{A}^*\bar{A} = \Fc^{-1} (\frac{\psi(\xi)}{\sqrt{e}})^2  \Fc $ and that
$$
f_{p/2\tau} \left(( \psi(\xi)/\sqrt{e})^2 \right) = \left( 1 + 2 |2\pi\xi|^{2 \tau} \int_0^T  \gamma(\lambda) \ud \lambda \right)^{-p/2 \tau},
$$
together with \eqref{eq5} yields \eqref{log sour cond}.\\
$\Leftarrow$) Conversely, assume \eqref{log sour cond} with $ q= p /2\tau$ then we get  
\begin{eqnarray*}
\label{eq6}
\hat{u}(\xi,0) & =& \left( 1 + 2 |2\pi\xi|^{2 \tau} \int_0^T  \gamma(\lambda) \ud \lambda \right)^{-p/2 \tau} \hat{w}(\xi) \\
& \leq & \min \left\lbrace 1, \left( 2 (2\pi)^{2 \tau} \int_0^T  \gamma(\lambda) \ud \lambda \right)^{-p/2\tau} \right\rbrace \left(1 + |\xi|^{2\tau} \right)^{-p/2\tau} |\hat{w}(\xi)|
\end{eqnarray*}
from which we deduce that $u(\cdot,0) \in H^p(\Rb^n)$ with $ ||u(\cdot,0)||_{H^p} \leq \tilde{C}(p) ||w||_{L^2}$ with $$\tilde{C}(p) = \min \left\lbrace 1, \left( 2 (2\pi)^{2 \tau} \int_0^T  \gamma(\lambda) \ud \lambda \right)^{-p/2\tau} \right\rbrace.$$
\end{Proof}

\begin{Remark}
From Proposition \ref{Prop 1 smoothness condition}, we can deduce that the smoothness condition \eqref{smoothness cond on u(cdot,0)} is nothing but the logarithmic source condition \eqref{log sour cond} whose order-optimal convergence rate under the noise level condition \eqref{noise level cond on data} is nothing but 
$C \rho f_{p/2\tau}(\delta^2/\rho^2)$ with $C \geq 1$ independent of $E$ and $\delta$. For details, see, e.g. \cite{hohage2000regularization}.
\end{Remark}
Now let us state the following proposition which gives a rate of convergence of the regularization error $u(\cdot ,0) - u_\beta$ under the smoothness condition \eqref{smoothness cond on u(cdot,0)}.
\begin{Proposition}
\label{Prop est rate reg error}
Consider the setting of Lemma \ref{lemma Alibaud et al}. Assume that the unknown solution $u(\cdot,0)$ satisfies the smoothness condition \eqref{smoothness cond on u(cdot,0)} and let $u_\beta$ be given in \eqref{def reg solution}. Then 
\begin{equation}
\label{estimate regul error}
|| u(\cdot,0) - u_\beta ||_{L^2}\leq \left(\frac{2s}{p\wedge 2s}\right)^q\rho f_q \left(\beta^{2s} \right) (1 + o(1)), \quad \text{as} \quad \beta \downarrow 0 \quad \text{with} \quad  q = p/2\tau,
\end{equation}
where $\rho = C(p) E$ in \eqref{log sour cond}.
\end{Proposition}
\begin{Proof}
Let $D_\beta = I-C_\beta$. From the first order optimality condition in \eqref{def reg solution}, and the fact that $A^*A = e\bar{A}^*\bar{A}$,  we get
\begin{eqnarray}
\label{eq7}
u(\cdot,0) - u_\beta & = &u(\cdot,0) - \left[ A^*A + D_\beta^*D_\beta \right]^{-1} A^* g \nonumber \\ 
& = &u(\cdot,0) -  \left[ A^*A + D_\beta^*D_\beta \right]^{-1} A^*A u(\cdot,0) \nonumber\\
& = &u(\cdot,0) -  \left[ \bar{A}^*\bar{A} + e^{-1}D_\beta^*D_\beta \right]^{-1} \bar{A}^*\bar{A} u(\cdot,0) \nonumber\\
& = &  e^{-1} \left[ \bar{A}^*\bar{A} + e^{-1}D_\beta^*D_\beta \right]^{-1} D_\beta^*D_\beta  u(\cdot,0).
\end{eqnarray}
From \eqref{eq7} and Proposition \ref{Prop 1 smoothness condition}, we get
\begin{equation}
\label{eq8}
u(\cdot,0) - u_\beta = e^{-1} \left[ \bar{A}^*\bar{A} + e^{-1}D_\beta^*D_\beta \right]^{-1} D_\beta^*D_\beta  f_q(\bar{A}^*\bar{A}) w.
\end{equation}
Now by noticing that operator $\bar{A}$ can be rewritten as a convolution operator (of kernel $\Fc^{-1}(\psi(\xi)/\sqrt{e})$), we get that the operator $C_\beta$ commutes with $\bar{A}^*\bar{A}$. By applying this commutation in \eqref{eq8}, we deduce that
\begin{equation}
\label{eq9}
u(\cdot,0) - u_\beta =  f_q(\bar{A}^*\bar{A})\, h_\beta,
\end{equation}
with 
\begin{equation}
\label{def h beta}
 h_\beta = e^{-1}\left[ \bar{A}^*\bar{A} + e^{-1}D_\beta^*D_\beta \right]^{-1}D_\beta^*D_\beta w.
\end{equation}
Using Parseval identity, we get that $||h_\beta||_{L^2} \leq ||w|| \leq  \rho $.
Hence, we have
\begin{equation}
\label{key part}
u(\cdot,0) - u_\beta =  f_q(\bar{A}^*\bar{A})\, h_\beta,\quad \text{with}\quad ||h_\beta||_{L^2}\leq  \rho.
\end{equation}
By application of \cite[Proposition 1]{hohage2000regularization} to the special case of logarithmic source function, and given that from \eqref{consistency result}, $|| T(u(\cdot,0) - u_\beta)|| \to 0$ as $\beta \downarrow 0$, we deduce that
\begin{equation}
\label{main1}
||u(\cdot,0) - u_\beta||_{L^2} \leq \rho f_q \left( \frac{|| \bar{A}(u(\cdot,0) - u_\beta)||_{L^2}^2}{\rho^2} \right) ( 1+ o(1)) \quad \text{as} \quad \beta \downarrow 0.
\end{equation} 
Using interpolation inequality, we have
\begin{eqnarray}
\label{eq9bis}
|| \bar{A}(u(\cdot,0) - u_\beta)|| & = & || (\bar{A}^*\bar{A})^{1/2}(u(\cdot,0) - u_\beta)||  \nonumber\\
& \leq & || (\bar{A}^*\bar{A})(u(\cdot,0) - u_\beta)||^{1/2} || u(\cdot,0) - u_\beta||^{1/2} .
\end{eqnarray}
On the one hand,
\begin{eqnarray}
\label{eq10}
|| (\bar{A}^*\bar{A})(u(\cdot,0) - u_\beta)|| & =& || (\bar{A}^*\bar{A})\left[ \bar{A}^*\bar{A} + e^{-1}D_\beta^*D_\beta \right]^{-1} e^{-1}D_\beta^*D_\beta  u(\cdot,0)|| \nonumber \\
& \leq & e^{-1} ||D_\beta^*D_\beta  u(\cdot,0)|| \quad \text{using Parseval identity} \nonumber \\
& \leq & e^{-1} \sqrt{C_2} \beta^{2(s\wedge(p/2)} E \quad \text{using} \,\, \eqref{upp bound (I-C_beta)^2 f} \,\, \text{and} \,\,\eqref{smoothness cond on u(cdot,0)} .
\end{eqnarray}
On the other hand, since $|| u(\cdot,0) - u_\beta ||  \to 0$ as $\beta \downarrow 0$, we have that for $\beta \ll 1$
\begin{equation}
\label{eq11}
|| u(\cdot,0) - u_\beta || \leq \frac{\rho^2}{e^{-1} \sqrt{C_2}  E}.
\end{equation} 
Putting together \eqref{main1}, \eqref{eq9bis} \eqref{eq10} and \eqref{eq11} and applying \eqref{Property log sourc function} yields \eqref{estimate regul error}.
\end{Proof}
\begin{Remark}
The commutation between the mollifier operator $C_\beta$ and operator $\bar{A}^*\bar{A}$ is crucial in the above proof. Indeed, without this commutation,we could not get the key estimate \eqref{eq9} from \eqref{eq8}.
\end{Remark}

Now we are ready to state on of the following theorem about order-optimality of the regularization method under the smoothness condition \eqref{smoothness cond on u(cdot,0)}.
\begin{Theorem}
\label{Theorem order optim 1 under noisy data only}
Consider the setting of Lemma \ref{lemma Alibaud et al}. Let $g^\delta \in L^2(\Rb^n)$ be a noisy distribution of final temperature satisfying  \eqref{noise level cond on data}. Assume that the solution $u(\cdot,0)$ satisfies \eqref{smoothness cond on u(cdot,0)} and let $u_\beta^\delta$ be the reconstructed solution defined by \eqref{def reg sol noisy data} using the noisy data $g^\delta$. Then for the a-priori selection rule $\beta(\delta) = \left( \Theta_q^{-1}(\delta/\rho) \right)^{1/2s}$ with $\Theta_q(t) = \sqrt{t}f_q(t)$ 
we have
\begin{equation}
\label{order optim rata noisy data only}
||u(\cdot,0) - u_{\beta(\delta)}^\delta||_{L^2} \leq K \rho f_q\left( \delta^2/\rho ^2 \right)(1 + o(1)) \quad \text{as} \quad \delta \to 0, 
\end{equation}
where $K$ is a constant independent of $E$ and $\delta$, $q = p/2\tau$ and $\rho = C(p)E$ given in \eqref{log sour cond}.
\end{Theorem}
\begin{Proof}
From Propositions \ref{Prop estimate data propagated error} and \ref{Prop 1 smoothness condition}, we deduce using triangular inequality that
$$
||u(\cdot,0) - u_{\beta}^\delta||_{L^2} \leq  \left(\frac{2s}{p\wedge 2s}\right)^q\rho f_q \left(\beta^{2s} \right) (1 + o(1)) + C \frac{\delta}{\beta^s} \quad \text{as} \quad \beta \downarrow 0.
$$
For $\beta(\delta) = \left( \Theta_q^{-1}(\delta/\rho) \right)^{1/2s} $, we deduce that
\begin{equation}
\label{eq12}
||u(\cdot,0) - u_{\beta(\delta)}^\delta||_{L^2} \leq  \left( \bar{C }\left(\frac{2s}{p\wedge 2s}\right)^q + C \right) \rho f_q \left(\Theta_p^{-1}(\delta/\rho) \right) 
\end{equation}
Estimate \eqref{order optim rata noisy data only} follows readily from the fact that $\rho f_q( \Theta_q^{-1}(\delta/\rho) )$ is nothing but the optimal rate under \eqref{log sour cond} (see, e.g. \cite[Theorem 1]{mathe2003geometry} and \cite[Theorem 2.1]{tautenhahn1998optimality}).
\end{Proof}

By the way, from Proposition \ref{Prop estimate data propagated error} and \ref{Prop est rate reg error}, we can readily establish the following result which exhibits a-priori parameter choice rule independent of the smoothness a-priori on the unknown solution.
\begin{Corollary}
Consider the setting of Theorem \ref{Theorem order optim 1 under noisy data only}. Then for the a-priori selection rule $\beta(\delta) = c \delta^{\theta/s}$ with $c>0$ and $\theta \in (0,1)$, we have
\begin{equation}
\label{order optim rata noisy data only}
||u(\cdot,0) - u_{\beta(\delta)}^\delta||_{L^2} =  \bigO{f_q(\delta^2)}, \quad \text{as} \quad \delta \to 0.
\end{equation}
\end{Corollary}
\begin{Remark}
We can see that the convergence rate in \eqref{order optim rata noisy data only} is actually order-optimal under the logarithmic source condition \eqref{log sour cond}. Hence Theorem \ref{Theorem order optim 1 under noisy data only} implies that the regularization method is order-optimal under the smoothness condition \eqref{smoothness cond on u(cdot,0)}. Moreover the a-priori selection rule given (i.e. $\beta(\delta) = c  \delta^{\theta/s}$ with $c>0$ and $\theta \in (0,1)$) is order-optimal with respect to the noise level $\delta$ and independent of the smoothness assumption on the solution $u(\cdot,0)$.
\end{Remark}

Now let us study error estimates when both the data and the operator are noisy.
\begin{Remark} Though the operator $A$ is explicitly known as \eqref{def operator A}, in practical implementation, this operator is approximated. For instance, given that integrals on unbounded domain are usually truncated numerically to a sufficiently large yet bounded domain, then in numerical implementation, the function $\psi$ which defines operator $A$, though having unbounded support is truncated on a bounded interval $I$. That is, the function $\psi$ is approximated by $1_I \psi$, so that $A$ is approximated by $\Fc^{-1} 1_I(\xi) \psi(\xi) \Fc$ where $1_I$ is the function equal $1$ on $I$ and $0$ outside of $I$.
\end{Remark}
 Henceforth, we set $A_h: L^2(\Rb^n) \to L^2(\Rb^n) $ to be a convolution operator approximating $A$ such that
\begin{equation}
\label{noisy cond on operator}
||| A - A_h||| \leq h, \quad \text{with} \quad |||A_h|||\leq 1.
\end{equation}
The requirement that $A_h$ is a convolution is quite reasonable in order to preserve the intrinsic property of the exact operator $A$. In the sequel, we call by $\psi_h$ the Fourier transform of the convolution kernel of operator $A_h$ so that $A_h = \Fc^{-1} \psi_h(\xi) \Fc$.

Let $u_\beta^{\delta,h}$ be the final reconstruction corresponding to the noisy data $g^\delta$ and noisy operator $A_h$ defined by
\begin{equation}
\label{def reg sol noisy data and noisy ope}
u_\beta^{\delta,h} = \mathrm{argmin}_{u \in L^2(\Rb^n)} || A_h u -g^\delta ||_{L^2}^2 + || (I-C_\beta)u||_{L^2}^2.
\end{equation}

\begin{Theorem}
\label{Theorem 2 order opti noisy data and ope}
Consider the setting of Lemma \ref{lemma Alibaud et al}. Let $g^\delta$ and $A_h$ satisfying \eqref{noise level cond on data} and \eqref{noisy cond on operator} respectively. Let the approximate solution $u_\beta^{\delta,h}$ be defined in \eqref{def reg sol noisy data and noisy ope}. Assume that the unknwon solution $u(\cdot,0)$ satisfies the smoothness condition \eqref{smoothness cond on u(cdot,0)}. Then for $\beta(\delta,h) = (h + \delta/\rho)^{1/2s}$ the following estimate holds
\begin{equation}
\label{order opt res under noisy data and oper}
|| u(\cdot,0) - u_{\beta(\delta,h)}^{\delta,h}||_{L^2} \leq \bar{K} \rho f_q\left( \left( h + \delta/\rho \right)^2 \right)
\end{equation}
where $q = p/2\tau$, $\rho = C(p)E$ given in \eqref{log sour cond} and $\bar{K}$ is a constant independent of $\rho$, $\delta$ and $h$.
\end{Theorem}
\begin{Proof}
For the sake of simplicity of notation, given a linear mapping $L: L^2(\Rb^n) \to L^2(\Rb^n)$, let us introduce the notation $S_\beta(L)$ and $R_\beta(L)$ for the linear mappings on $L^2(\Rb^n)$ defined respectively by
\begin{equation}
\label{notation s beta}
\forall f \in L^2(\Rb^n), \quad S_\beta(L)f : =  \left[ L^*L + D_\beta^*D_\beta \right]^{-1} L^* f,
\end{equation}
and
\begin{equation}
\label{notation R beta}
\forall f \in L^2(\Rb^n), \quad R_\beta(L)f : =  \left[ L^*L + D_\beta^*D_\beta) \right]^{-1} D_\beta^*D_\beta f,
\end{equation}
where $D_\beta = I - C_\beta$.
By Applying the first order optimality condition in \eqref{def reg sol noisy data and noisy ope}, we get $u_{\beta}^{\delta,h}  = S_\beta(A_h) g^\delta $. Let $w \in L^2(\Rb^n)$ be the function given in \eqref{log sour cond} and $\bar{A_h} = A_h/\sqrt{e}$ be a normalized version of $A_h$ such that $|||\bar{A_h}^*\bar{A_h}|||\leq \exp(-1)$, we have
\begin{eqnarray}
\label{eq 13}
|| u(\cdot,0) - u_{\beta}^{\delta,h} || & \leq & || u(\cdot,0) - S_\beta(A_h) A_h u(\cdot,0) || + || S_\beta(A_h) A_h u(\cdot,0)  - S_\beta(A_h) g^\delta|| \nonumber \\
& = & || R_\beta(A_h) u(\cdot,0) || + || S_\beta(A_h) (A_h u(\cdot,0)  - g^\delta)|| \nonumber \\
& \leq & || R_\beta(A_h) f_q(\bar{A_h}^*\bar{A_h})w || + || R_\beta(A_h) (u(\cdot,0) - f_q(\bar{A_h}^*\bar{A_h})w) || \nonumber\\
&  & \quad + \quad || S_\beta(A_h) (A_h u(\cdot,0)  - g^\delta)|| \nonumber \\
& \leq & || R_\beta(A_h) f_q(\bar{A_h}^*\bar{A_h})w || + || R_\beta(A_h) ( f_q(\bar{A}^*\bar{A}) - f_q(\bar{A_h}^*\bar{A_h}))(w) || \\
&  & \quad + \quad || S_\beta(A_h) (A_h u(\cdot,0)  - g^\delta)||.  \nonumber
\end{eqnarray}
But
\begin{equation}
\label{eq15}
R_\beta(A_h) f_q(\bar{A_h}^*\bar{A_h})w  = e^{-1}\left[ \bar{A_h}^*\bar{A_h} + e^{-1}D_\beta^*D_\beta \right]^{-1}D_\beta^*D_\beta f_q(\bar{A_h}^*\bar{A_h}) w \quad \text{with} \quad ||w || \leq \rho.
\end{equation}
Given that $A_h$ is a convolution operator, then $A_h^*A_h$ commutes with $C_\beta$, then \eqref{eq15} yields
\begin{equation}
\label{eq16}
R_\beta(A_h) f_q(\bar{A_h}^*\bar{A_h})w =  f_q(\bar{A_h}^*\bar{A_h}) \bar{h_\beta} \quad \text{with} \quad ||\bar{h_\beta} || \leq ||w|| \leq  \rho, \quad q= p/2\tau.
\end{equation}
More precisely $\bar{h_\beta} = e^{-1}\left[ \bar{A_h}^*\bar{A_h} + e^{-1}D_\beta^*D_\beta \right]^{-1}D_\beta^*D_\beta w$.
Following the same lines as the proof of \eqref{estimate regul error} at the only difference that $A$ is replaced by $A_h$, we get that
\begin{equation}
\label{eq17}
|| R_\beta(A_h) f_q(\bar{A_h}^*\bar{A_h})w || \leq \left(\frac{2s}{p\wedge 2s}\right)^q\rho f_q \left(\beta^{2s} \right)  (1 + o(1)) \quad \text{as} \quad \beta \downarrow 0.
\end{equation}
Next, using \cite[Lemma 9]{hohage2000regularization}, we get that
\begin{equation}
||| f_q(\bar{A}^*\bar{A}) - f_q(\bar{A_h}^*\bar{A_h}) ||| \leq 2 f_q \left(||| \bar{A}^*\bar{A} - \bar{A_h}^*\bar{A_h}||| \right).
\end{equation}
But using \eqref{noisy cond on operator},
\begin{equation}
\label{eq18}
||| \bar{A}^*\bar{A} - \bar{A_h}^*\bar{A_h} ||| \leq ||| \bar{A}^* ( \bar{A} - \bar{A_h}) ||| + ||| (\bar{A}^* - \bar{A_h}^*)\bar{A_h}||| \leq \frac{2}{e} ||| A_h - A || \leq \frac{2}{e} h.
\end{equation}
Thus \eqref{eq17}, \eqref{eq18} and the fact that $|||R_\beta(A_h)|||\leq 1$ yields
\begin{equation}
\label{eq19}
|| R_\beta(A_h) ( f_q(\bar{A}^*\bar{A}) - f_q(\bar{A_h}^*\bar{A_h}))(w) || \leq 2 \rho f_q\left( \frac{2 h}{e} \right) \leq 2 \rho f_q(h).
\end{equation}
Finally, following the same lines as equations \eqref{eq3} to \eqref{zzz mol deconv} except that the function $\psi$ (resp. $g - g^\delta$) is replaced by $\psi_h$ (resp. $A_h u(\cdot,0)  - g^\delta$), one gets
\begin{eqnarray}
\label{eq20}
|| S_\beta(A_h) (A_h u(\cdot,0)  - g^\delta)|| & \leq & C \frac{||A_h u(\cdot,0)  - g^\delta||}{\beta^s} \nonumber \\
& \leq & C \frac{||A_h u(\cdot,0)- A u(\cdot,0)|| + ||A u(\cdot,0) - g^\delta||}{\beta^s}  \nonumber \\
& \leq & C\,  \frac{h \rho + \delta}{\beta^s} \quad \text{using}\quad \eqref{noisy cond on operator}, \eqref{log sour cond}, \eqref{noise level cond on data}.
\end{eqnarray}
Thus, from \eqref{eq 13}, \eqref{eq17}, \eqref{eq19} and \eqref{eq20}, we get
\begin{equation}
\label{eq21}
|| u(\cdot,0) - u_{\beta}^{\delta,h} || \leq  \left(\frac{2s}{p\wedge 2s}\right)^q \rho f_q \left(\beta^{2s} \right) (1 + o(1)) + 2 \rho f_q(h) + C \,\frac{h \rho + \delta}{\beta^s} \quad \text{as} \quad \beta \downarrow 0.
\end{equation}
Let $\beta : = (h + \delta/\rho)^{1/2s}$, and $h,\rho, \delta$ such that $h + \delta/\rho \ll 1$, applying Lemma \ref{Lemma prop log sour func}, we get
\begin{equation}
\label{eq22}
\begin{cases}
\vspace{0.4cm}
\rho f_q \left(\beta^{2s} \right) = \rho f_q \left(h + \delta/\rho\right)\leq 2^q \rho f_q \left( (h + \delta/\rho)^2 \right) \\
\vspace{0.4cm}
\rho f_q(h) \leq 2^q \rho f_q(h^2) \leq 2^q \rho f_q \left( (h + \delta/\rho)^2 \right) \\
\vspace{0.4cm}
C \,\frac{h \rho + \delta}{\beta^s} = C \,\rho \sqrt{h + \delta/\rho} = \rho\, o\left(f_q \left( (h + \delta/\rho)^2 \right) \right) \quad \text{as} \quad h + \delta/\rho \to 0.
\end{cases}
\end{equation}
Estimate \eqref{order opt res under noisy data and oper} follows readily from \eqref{eq21} and \eqref{eq22}.
\end{Proof}

Now, let us turn to the practical aspect of the choice of the regularization parameter $\beta$.
\section{Parameter choice rule}\label{section par choice rule}
In the implementation of a regularization method, a very important step is the choice of the regularization parameter. Usually, a-posteriori parameter choice rules, i.e. parameter choices depending both on the noisy data $g^\delta$ and the noise level $\delta$ are advocated.
In this section we present an order-optimal a-posteriori parameter choice rule closely related to Morozov principle. 

Given $r \in (0,1]$, a noisy data $g^\delta$ satisfying \eqref{noise level cond on data} and the approximate solution $u_\beta^\delta$ defined in \eqref{def reg sol noisy data}, let $\beta(\delta,g^\delta,r)$ be expressed as
\begin{equation}
\label{def beta a posteriori rule}
\beta(\delta,g^\delta,r) = \sup \left\lbrace \beta > 0, \quad \text{s.t.} \quad || A u_\beta^\delta - g^\delta|| < \delta + \delta^r  \right\rbrace.
\end{equation} 
Before moving to error estimates, let us first discuss the existence and uniqueness of $\beta(\delta,g^\delta,r)$ defined in \eqref{def beta a posteriori rule}.
Let us state the following result of existence and uniqueness of the parameter $\beta(\delta,g^\delta,r)$ defined in \eqref{def beta a posteriori rule}.
\begin{Proposition}
\label{Prop exist post par choice rule}
Assume that the noise level $\delta$ and the noisy data $g^\delta$ satisfies
\begin{equation}
\label{noise ration cond}
\delta + \delta^r \leq \frac{1}{2} ||g^\delta||,
\end{equation}
then the parameter $\beta(\delta,g^\delta,r)$ expressed in \eqref{def beta a posteriori rule} is well defined and satisfies
\begin{equation}
\label{char beta a posteriori rule}
|| A u_{\beta(\delta,g^\delta,r)}^\delta - g^\delta|| = \delta + \delta^r .
\end{equation}
\end{Proposition}
\begin{Proof}
Let $g^\delta$ be a noisy data satisfying \eqref{noise level cond on data} and \eqref{noise ration cond}. Let us introduce the function $\L(\beta): = ||A u_\beta^\delta - g^\delta||^2 $. Using Parseval identity, we have
\begin{equation}
\label{def func l beta}
\L(\beta) =  \int_{\Rb^n} |\Pi(\beta,\xi)|^2 |\widehat{g^\delta}(\xi)|^2 \ud \xi, \quad \text{with} \quad \Pi(\beta,\xi) = \frac{| 1- \widehat{\phi}(\beta\xi)|^2}{|\psi(\xi)|^2 + | 1- \widehat{\phi}(\beta\xi)|^2}.
\end{equation}
Given that for all $\xi \in \Rb^n$ and for all $\beta \in \Rb^+$, $\Pi(\beta,\xi) \leq 1$ and that $\forall \xi \in \Rb^n,\,\,\,\Pi(\beta,\xi) \to 0 $ as $\beta \downarrow 0$, we deduce that $\L(\beta) \to 0$ as $\beta \downarrow 0$ using the dominated convergence theorem. Let $\beta_1, \beta_2 \in \Rb_+^*$, by mere computation, one gets
\begin{equation}
\label{efg}
\L(\beta_1) - \L(\beta_2) = \int_{\Rb^n}  X(\beta_1,\beta_2,\xi) \left(| 1- \widehat{\phi}(\beta_1\xi)|^2 -| 1-\widehat{\phi}(\beta_2\xi)|^2 \right)
|\widehat{g^\delta}(\xi)|^2 \ud \xi,
\end{equation}
where
$$
X(\beta_1,\beta_2,\xi) = \frac{| \Pi(\beta_1,\xi) + \Pi(\beta_2,\xi)|\,\psi(\xi)}{|\psi(\xi)^2 + | 1- \widehat{\phi}(\beta_1\xi)|^2 | \times|\psi(\xi)^2 + | 1- \widehat{\phi}(\beta_2\xi)|^2 |}.
$$
Since $\widehat{\phi}$ is radially decreasing, then \eqref{efg} implies for $\beta_1 > \beta_2$, $\L(\beta_1) - \L(\beta_2)>0$ which shows that the function $\L$ is strictly increasing. Finally using Fatou Lemma and the fact that $||\psi||_\infty \leq 1$, one gets
\begin{eqnarray*}
\label{eq25}
\lim_{\beta \to + \infty} \L(\beta) \geq  \int_{\Rb^n} \lim_{\beta \to +\infty} \Pi(\beta,\xi)^2 | \widehat{g^\delta}(\xi)|^2 \ud \xi
= \int_{\Rb^n} \left( \frac{1}{|\psi(\xi)|^2 +1 } \right)^2 | \widehat{g^\delta}|^2 \ud \xi
 \geq  \int_{\Rb^n} \frac{1}{4} | \widehat{g^\delta}|^2 \ud \xi.
\end{eqnarray*}
Hence, in summary, the function $\L$ is continuous, strictly increasing and satisfies $ \lim_{\beta \to 0} \L(\beta) = 0$ and $\lim_{\beta \to \infty}\L(\beta) > ||g^\delta||^2/4$. This proves that under noise ratio condition \eqref{noise ration cond}, there exists a unique $\beta >0$ satisfying \eqref{def beta a posteriori rule} and that such $\beta$ is characterized by \eqref{char beta a posteriori rule}.
\end{Proof}

Notice that the condition \eqref{noise ration cond} is merely saying that the noisy data $g^\delta$ is not dominated by noise, for otherwise it is hopeless to recover meaningful approximate solution. 

\begin{Lemma}
\label{Lemma rates under opsterioru rule}
Consider the setting of Lemma  \ref{lemma Alibaud et al}. Let the solution $u(\cdot,0)$ of equation \eqref{op eq of our pb} satisfies the smoothness condition \eqref{smoothness cond on u(cdot,0)}. Let $g^\delta$ be a noisy data verifying \eqref{noise level cond on data} and \eqref{noise ration cond} and $u_\beta$ and $u_\beta^\delta$ be the regularized solutions defined by \eqref{def reg solution} and \eqref{def reg sol noisy data} respectively. Let $\rho = C(p)E $ be given in \eqref{log sour cond}, $q = p/2\tau$ and $\beta_r: = \beta(\delta,g^\delta,r)$ be defined in \eqref{def beta a posteriori rule}, then there exist constants $K_1(p)$ and $K_2(p)$ independent of $\delta$ and $\rho$ such that
\begin{equation}
\label{est reg err posteriori}
|| u(\cdot,0) - u_{\beta_r}||_{L^2} \leq   K_1(p)\, \rho f_q(\delta^{2}/\rho^2)(1 + o(1)) \quad \text{as} \quad \delta \to 0.
\end{equation}
and 
\begin{equation}
\label{est data propag err posteriori}
|| u_{\beta_r} - u_{\beta_r}^\delta ||_{L^2} \leq K_2(p) \rho^{\frac{s}{s \wedge (p/2)}} \delta^{1 - \frac{r\,s}{s \wedge (p/2)}} \left( f_{q} \left(\delta^2/\rho^2 \right) \right)^{\frac{s}{2(s \wedge (p/2))}} (1 + o(1)) \,\, \text{as} \,\, \delta \to 0
\end{equation}
\end{Lemma}
\begin{Proof} Let $\beta_r := \beta(\delta,g^\delta,r)$ defined in \eqref{def beta a posteriori rule} and $u(\cdot,0)$ satisfying \eqref{smoothness cond on u(cdot,0)}. From \eqref{key part}, we know that $u(\cdot,0) - u_{\beta_r} = f_q(\bar{A}^*\bar{A})\bar{h_{\beta_r}}$ with $||\bar{h_{\beta_r}}|| \leq \rho$. Let us consider the notation $S_\beta$ and $R_\beta$ defined in \eqref{notation s beta} and \eqref{notation R beta}. We have
\begin{eqnarray*}
|| A (u(\cdot,0) - u_{\beta_r}) || & = & || R_{\beta_r}(A^*) A u(\cdot,0) || \\
& \leq & || R_{\beta_r}(A^*) (A u(\cdot,0) - g^\delta) || +  || R_{\beta_r}(A^*) g^\delta|| \\
& =& || R_{\beta_r}(A^*) (g - g^\delta) || +  || A u_{\beta_r}^\delta - g^\delta|| \\
& \leq & \delta + \delta + \delta^r \quad \text{using that} \quad ||| R_{\beta_r}(A^*)||| \leq 1,\,\, \eqref{noise level cond on data} \,\, \text{and} \,\, \eqref{char beta a posteriori rule}.
\end{eqnarray*}
Then, we have 
\begin{equation*}
\begin{cases}
u(\cdot,0) - u_{\beta_r} = f_q(\bar{A}^*\bar{A})\bar{h_{\beta_r}} & \text{with} \,\, ||\bar{h_{\beta_r}}|| \leq \rho \\
|| \bar{A} (u(\cdot,0) - u_{\beta_r}) || \leq \delta^r (1 + 2\delta^{1-r})/\sqrt{e},
\end{cases}
\end{equation*}
which implies that 
\begin{equation}
\label{eq31}
|| u(\cdot,0) - u_{\beta_r}|| \leq  \omega(\delta^r(1 + 2\delta^{1-r})/\sqrt{e}, M_{f_q}(\rho),\bar{A})
\end{equation}
where for $\epsilon>0$, $\omega(\epsilon, M_{f_q}(\rho),\bar{A}) = \sup \left\lbrace ||f|| , \,\,\, f \in M_{f_q}(\rho), \,\,\, \text{s.t.} \,\,\, || \bar{A}f || \leq \epsilon \right\rbrace $ denotes the modulus of continuity of operator $\bar{A}$ on the subspace 
$$
M_{f_q}(\rho) = \left\lbrace f_q(\bar{A}^*\bar{A})w, \,\,\, w \in L^2(\Rb^n), \,\,\,\text{with} \,\,\,||w|| \leq \rho \right\rbrace.
$$
Given that $\omega(\epsilon, M_{f_q}(\rho),\bar{A}) \leq \rho f_q(\epsilon^{2}/\rho
^2)(1 + o(1))$ as $\epsilon \to 0$ (see, e.g. \cite[Proposition 2]{hohage2000regularization}) and that $\omega(\alpha \epsilon, M_{f_q}(\rho),\bar{A}) \leq \alpha \omega(\epsilon, M_{f_q}(\rho),\bar{A})$, from \eqref{eq31}, we deduce that
\begin{equation}
\label{eq32}
|| u(\cdot) - u_{\beta_r}|| \leq  ((1 + 2\delta^{1-r})/\sqrt{e}) \rho f_q(\delta^{2r}/\rho
^2)(1 + o(1)) \leq 2 \rho f_q(\delta^{2r}/\rho
^2)(1 + o(1)) \quad \text{as} \quad \delta \to 0.
\end{equation}
Applying \eqref{Property log sourc function} with $a=2r$, $b = 2$, $\lambda = 1/\rho^2$ and $t = \delta$ to \eqref{eq32} yields \eqref{est reg err posteriori}.

Now let us find a lower bound of $\beta_r$ in function of $\delta$ in order to deduce the rate in the propagated data noise error using \eqref{bound error prog term}. Let $\bar{\beta} = \lambda \beta_r$  with $\lambda >1$. We have
\begin{eqnarray}
\label{eq33}
|| A (u(\cdot,0) - u_{\bar{\beta}}) || & \geq & || A u_{\bar{\beta}}^\delta - g^\delta || - || A (u_{\bar{\beta}}^\delta - u_{\bar{\beta}}) - ( g^\delta - g) || \nonumber \\
& = & || A u_{\bar{\beta}}^\delta - g^\delta || - || R_{\bar{\beta}}(A^*) ( g^\delta - g) || \nonumber \\
& \geq & \delta + \delta^r - \delta = \delta^r \quad \text{from}  \quad\eqref{def beta a posteriori rule} \quad \text{and the fact that} \quad \bar{\beta} > \beta_r.
\end{eqnarray}
On the other hand, from \eqref{log sour cond} and the interpolation inequality, we have
\begin{equation}
\label{eq34}
|| A (u(\cdot,0) - u_{\bar{\beta}}) || = || (A^*A)^{1/2} R_{\bar{\beta}}(A) u(\cdot,0)  || \leq || A^*A R_{\bar{\beta}}(A) u(\cdot,0)  ||^{1/2} ||R_{\bar{\beta}}(A) u(\cdot,0)||^{1/2}
\end{equation}
Using Parseval identity, \eqref{log sour cond} and \eqref{smoothness cond on u(cdot,0)}, we have
\begin{equation}
\label{eq35}
 || A^*A R_{\bar{\beta}}(A) u(\cdot,0)  || \leq || D_{\bar{\beta}}^*D_{\bar{\beta}} u(\cdot,0)  ||  \leq \sqrt{C_2}\bar{\beta}^{2(s \wedge (p/2))} E = (\sqrt{C_2}/C(p))\,\rho \,\bar{\beta}^{2(s \wedge (p/2))} .
\end{equation}
Moreover, using \eqref{estimate regul error}, we get
\begin{eqnarray}
\label{eq36}
||R_{\bar{\beta}}(A) u(\cdot,0)|| & = & || u(\cdot,0) - u_{\bar{\beta}}|| \nonumber \\
& \leq & K' \rho f_q\left( \bar{\beta}^{ s \wedge (p/2)}\right)(1 + o(1)) \quad \text{as} \quad \bar{\beta} \downarrow 0\,\,\, 
\text{using} \,\,\, \eqref{Property log sourc function}
\end{eqnarray}
with $K'$ independent of $\rho$ and $\bar{\beta}$.
From \eqref{eq33} to \eqref{eq36}, we deduce that
\begin{equation}
\label{eq37}
\delta^r \leq || A (u(\cdot,0) - u_{\bar{\beta}}) || \leq  \bar{K}(p) \rho \, \Theta(\bar{\beta}^{s \wedge (p/2)}),
\end{equation}
where $\Theta$ is the function defined on $(0,1]$ by $\Theta(t) = t f_{q/2} ( t )$, 
and $K(p)$  is a constant independent of $\bar{\beta}$, $\rho$ and  $\delta$.
The function $\Theta$ is monotonically increasing, bijective from $(0,1]$ to $(0,1]$ with inverse given in \cite[Lemma 3.3]{tautenhahn1998optimality} as
$$
\Theta^{-1}(t) = t (- \ln t)^{q/2}(1 + o(1))  \quad \text{as} \quad t \to 0.
$$
Then, as $\delta \to 0$, applying the function $\Theta^{-1}$ to \eqref{eq37} yields
\begin{equation}
\label{eq38}
 \left( \frac{\delta^r}{\bar{K}(p) \rho} \right) \left( - \ln \left( \frac{\delta^r}{\bar{K}(p) \rho}\right) \right)^{q/2}(1 + o(1))  \leq ( \lambda {\beta_r})^{s \wedge (p/2)}, \quad \text{as} \quad \delta \to 0 .
\end{equation} 
Without loss of generality, assuming that $\bar{K}(p)\geq 1$, \eqref{eq38} implies that
\begin{equation}
\label{eq39}
 \left( \frac{\delta^r}{\bar{K}(p) \rho} \right) \left( - \ln ( \delta^r/\rho) \right)^{q/2}(1 + o(1))  \leq (\lambda {\beta_r})^{s \wedge (p/2)}, \quad \text{as} \quad \delta \to 0 .
\end{equation}
Hence from Proposition \ref{Prop estimate data propagated error} and \eqref{eq39} we deduce that
\begin{eqnarray}
\label{eq40}
|| u_{\beta_r} - u_{\beta_r}^\delta || & \leq & C  \frac{\delta}{(\beta_r)^s} \nonumber \\
& \leq & C \lambda^s \bar{K}(p)^{\frac{s}{s \wedge (p/2)}} \rho^{\frac{s}{s \wedge (p/2)}} \delta^{1 - \frac{r\,s}{s \wedge (p/2)}} \left( f_{q} \left(\delta^r/\rho \right) \right)^{\frac{s}{2(s \wedge (p/2))}} (1 + o(1)) \,\, \text{as} \,\, \delta \to 0 \nonumber \\
& \leq & K_2(p) \rho^{\frac{s}{s \wedge (p/2)}} \delta^{1 - \frac{r\,s}{s \wedge (p/2)}} \left( f_{q} \left(\delta^2/\rho^2 \right) \right)^{\frac{s}{2(s \wedge (p/2))}} (1 + o(1)) \,\, \text{as} \,\, \delta \to 0,
\end{eqnarray}
where the last inequality uses \eqref{Property log sourc function 2} and \eqref{Property log sourc function}.
\end{Proof}

The following theorem exhibits the convergence rates obtained by the a-posteriori rule \eqref{def beta a posteriori rule} under the smoothness condition \eqref{smoothness cond on u(cdot,0)}.
\begin{Theorem}
\label{Theorem conv rate for a posteriori rule}
Consider the setting of Lemma \ref{Lemma rates under opsterioru rule} and let $\beta_r := \beta(\delta,g^\delta,r)$ given in \eqref{def beta a posteriori rule} with $r \in (0,1]$. Let $s$ be the parameter in \eqref{cond on conv kernel phi}. The following holds
\begin{itemize}
\item[i)] If $s \leq p/2$, then for all $r \in (0,1)$, the parameter selection rule \eqref{def beta a posteriori rule} is order-optimal under the smoothness condition \eqref{smoothness cond on u(cdot,0)}, that is, there exists a constant $K_3(p)$ independent of $\rho$ and $\delta$ such that
\begin{equation}
\label{conv rate for a posteriori rule}
|| u(\cdot,0) - u_{\beta_r}^\delta || \leq K_3(p) \rho f_q(\delta^2/\rho^2)(1 + o(1)) \quad \text{as} \quad \delta \to 0.
\end{equation}

\item[ii)] Assume that $s > p/2$. Then for all $r < p/2s$, there exists a constant $K_4(p)$ independent of $\rho$ and $\delta$ such that
\begin{equation}
\label{optimal rate for a posteriori rule up to power}
|| u(\cdot,0) - u_{\beta_r}^\delta || \leq K_4(p) \rho^{2s/p} f_q(\delta^2/\rho^2)(1 + o(1)) \quad \text{as} \quad \delta \to 0.
\end{equation}
For $ r = p/2s $ there exists a constant $K_5(p)$ independent of $\rho$ and $\delta$ such that
\begin{equation}
\label{sub optimal rate for a posteriori rule}
|| u(\cdot,0) - u_{\beta_r}^\delta || \leq K_5(p) \rho^{2s/p} f_q(\delta^2/\rho^2)^{1\wedge (s/p)}(1 + o(1)) \quad \text{as} \quad \delta \to 0.
\end{equation}
\end{itemize}
\end{Theorem}
\begin{Proof} From \eqref{est reg err posteriori}, we deduce the a-posteriori rule \eqref{def beta a posteriori rule} is order-optimality under \eqref{smoothness cond on u(cdot,0)} if the rate in \eqref{est data propag err posteriori} of the propagated data noise error is also order-optimal.\\
i) If $s \leq p/2$, and $r\in(0,1)$ then $\frac{s}{s \wedge (p/2)} = 1$ and $ \frac{s}{2(s\wedge (p/2))} = 1/2$  which implies that 
\begin{equation}
\label{eqs}
\rho^{\frac{s}{s \wedge (p/2)}} \delta^{1 - \frac{r\,s}{s \wedge (p/2)}} \left( f_{q} \left(\delta^2/\rho^2 \right) \right)^{\frac{s}{2(s \wedge (p/2))}} = \rho\, \delta^{1-r} \left( f_{q} \left(\delta^2/\rho^2 \right) \right)^{1/2} = \rho \, o\left( f_q(\delta^2/\rho^2) \right) \quad \text{as} \quad \delta \to 0.
\end{equation}
ii) Assume that $s > p/2$, then 
$$
\rho^{\frac{s}{s \wedge (p/2)}} \delta^{1 - \frac{r\,s}{s \wedge (p/2)}} \left( f_{q} \left(\delta^2/\rho^2 \right) \right)^{\frac{s}{2(s \wedge (p/2))}} = \rho^{\frac{2s}{p}} \delta^{1 - \frac{2 r \,s}{p}} \left( f_{q} \left(\delta^2/\rho^2 \right) \right)^{s/p}. 
$$
Therefore, if $r  < p/2s$ then $1 - \frac{2 r \,s}{p} >0$ which implies that
$$
\rho^{\frac{2s}{p}} \delta^{1 - \frac{2 r \,s}{p}} \left( f_{q} \left(\delta^2/\rho^2 \right) \right)^{s/p} = \rho^{\frac{2s}{p}} o\left( f_q(\delta^2/\rho^2) \right) \quad \text{as} \quad \delta \to 0,
$$
from which \eqref{optimal rate for a posteriori rule up to power} follows.
Finally if $r =p/2s$, then 
$$
\rho^{\frac{s}{s \wedge (p/2)}} \delta^{1 - \frac{r\,s}{s \wedge (p/2)}} \left( f_{q} \left(\delta^2/\rho^2 \right) \right)^{\frac{s}{2(s \wedge (p/2))}} = \rho^{\frac{2s}{p}}  \left( f_{q} \left(\delta^2/\rho^2 \right)  \right)^{s/p}  \quad \text{as} \quad \delta \to 0.
$$
from which \eqref{sub optimal rate for a posteriori rule} follows.
\end{Proof}

\begin{Remark}
Except the fact the power of $\rho$ in \eqref{optimal rate for a posteriori rule up to power} is $2s/p >1$, the rate given in \eqref{optimal rate for a posteriori rule up to power} may be qualified as `\it{order-optimal}'. Similarly if $s \geq p$ the rate in \eqref{sub optimal rate for a posteriori rule} may also be qualified as `\it{order-optimal}'.
\end{Remark}
\begin{Remark}{\label{remark 6}}
From Theorem \ref{Theorem conv rate for a posteriori rule}, we get that the a-posteriori rule \eqref{def beta a posteriori rule} is order-optimal for all $r \in (0,1)$ if $s \leq p/2$. From \eqref{eqs}, we can also deduce that if $s \leq p/2$ and $r=1$, then we get the sub-optimal rate $\bigO{\rho\, \left( f_{q} \left(\delta^2/\rho^2 \right) \right)^{1/2} }$ as $\delta \to 0$.
\end{Remark}

We end up this section by the following algorithm for estimating the parameter $\beta_r$. In Algorithm \ref{Algo beta r}, the function $\Pi$ is defined by \eqref{def func l beta}.
\begin{algorithm}
\begin{center}
\begin{algorithmic}[1]
\State Set $\beta_0 \gg 1$ and $q \in (0,1)$
\State Set $\beta_r = \beta_0$ (initial guess)
\While{$||\Pi(\beta_r,\xi) \widehat{g^\delta}(\xi)  ||_{L^2} \,>\, \delta + \delta^r$}
\State $\beta_r = q \times \beta_r$
\EndWhile
\end{algorithmic}
\end{center}
\caption{}
\label{Algo beta r}
\end{algorithm}

\section{Numerical experiments}\label{section numerical experiments}

 In order to illustrate the efficiency and robustness of the regularization method presented in this paper, we treat four numerical examples in two-dimension space where the final time $T$ is invariably set to $1$.

\textbf{Example 1}: $\gamma(t) = 0.1(3 - 2t)$ and $u(x,0) = e^{-x_1^2 - x_2^2}$.

\textbf{Example 2}: $\gamma(t) = 0.1(2 -t) $ and
$
u(x,0) =v(x_1)v(x_2)$ where $v$ is triangle impulse defined by $v(\lambda)=
\begin{cases}
1 + \lambda/3 & \text{if} \,\, \lambda \in [-3,0]\\
1 - \lambda/3 & \text{if} \,\, \lambda \in (0,3] \\
0 & \text{otherwise}.
\end{cases}
$

\textbf{Example 3}:  $\gamma(t) = 0.1 $ and
$
u(x,0) =
\begin{cases}
1  & \text{if} \,\, (x_1,x_2) \in [-5,5]^2\\
0 & \text{otherwise}.
\end{cases}
$

\textbf{Example 4}: we consider an image deblurring process where $u(\cdot,0)$ is the head phantom image. We set $\tau = 1/2$ which has already been used for modeling X-ray scattering. Here $\gamma$ is set to $0.1$. 

In Examples 1, 2 and 3, we consider classical diffusion, i.e. $\tau = 1$.

The initial temperature distribution $u(\cdot,0)$ in Examples 2 and 3 are two-dimensional versions of examples given in \cite{hao2011stability}. Notice that the smoothness of the solution is decreasing from Example 1 to Example 4. Indeed, in Example 1, $u(\cdot,0) \in H^{p}(\Rb^2)$ for all $p>0$; in Example 2, $u(\cdot,0) \in H^{1}(\Rb^2)$; in Example 3, $u(\cdot,0) \in H^{p}(\Rb^2)$ for $p <1/2$ while in Example 4, $u(\cdot,0)$ is very irregular and does not belongs to any Sobolev space $H^p(\Rb^p)$ with $p>0$. Therefore, under a same setting, the quality of reconstruction is expected to get better as we move from Example 3 to Example 1. 

In the four examples, the support of $u(\cdot,0)$ is $[-L,L]^2$ which is uniformly discretized as $(x_1(i),x_2(j))$ where $x_1(i)=x_2(i) = -L + (i-0.5)\kappa$ with $\kappa = 2L/N$ and $i,j=1,...,N$. In all the simulations, we set $L=10$ and $N =256$.

The noisy data $g^\delta$ is generated as $g^\delta(x_1(i),x_2(j)) = u(x_1(i),x_2(j),T) + \eta \epsilon(x_1(i),x_2(j))
$ where $\epsilon(x_1(i),x_2(j))$ is a random number drawn from the standard normal distribution. Given a noise level $\delta$, $\eta$ is set such that $\eta  E \left[ ||\epsilon||_2 \right] = \delta$. In the simulation, given a percentage of noise $perc\_noise$, we set $\delta = perc\_noise \times||u(\cdot,T)||_2/100$.

The mollifier operator $C_\beta$ uses the standard normal convolution kernel, that is,
$$
\phi(x) = (1/2\pi) \exp(-(x_1^2+x_2^2)/2),
$$
 which satisfies condition \eqref{cond on conv kernel phi} (with $s=2$).

Given the noisy data $g^\delta$, we computed its Fourier transform $\widehat{g^\delta}$ and the reconstructed solution $u_\beta^\delta$ is computed as the inverse Fourier transform of the function $\Lambda_\beta^\delta$ defined by
$$
\Lambda_\beta(\xi) =  \frac{\psi(\xi) \widehat{g^\delta}(\xi)}{\psi^2(\xi) + |1 - \hat{\phi}(\beta \xi)|^2}, \quad \xi \in \Rb^2.
$$
Notice that the function $\psi$ for examples 1 to 4 are given by 
$$
\psi_1(\xi) = e^{- 0.8 \pi^2 |\xi|^2 }, \,\, \psi_2(\xi) = e^{- 0.6\pi^2 |\xi|^2}, \,\, \psi_3(\xi) = e^{- 0.4\pi^2 |\xi|^2 } \,\, \text{and} \,\, \psi_4(\xi)= e^{- 0.2\pi |\xi| }.
$$
The Fourier transform and inverse Fourier transform involved in the computation of the reconstructed solution $u_\beta^\delta$ are quite rapidly evaluated with the numerical procedure from \cite{bailey1994fast} using fast Fourier transform (FFT) algorithm. From the Shannon-Nyquist principle, we set the frequency domain corresponding to $[-\Omega,\Omega]^2$ with $\Omega = N/4L$.

\begin{figure}[h!]
\begin{center}
\includegraphics[scale=0.4]{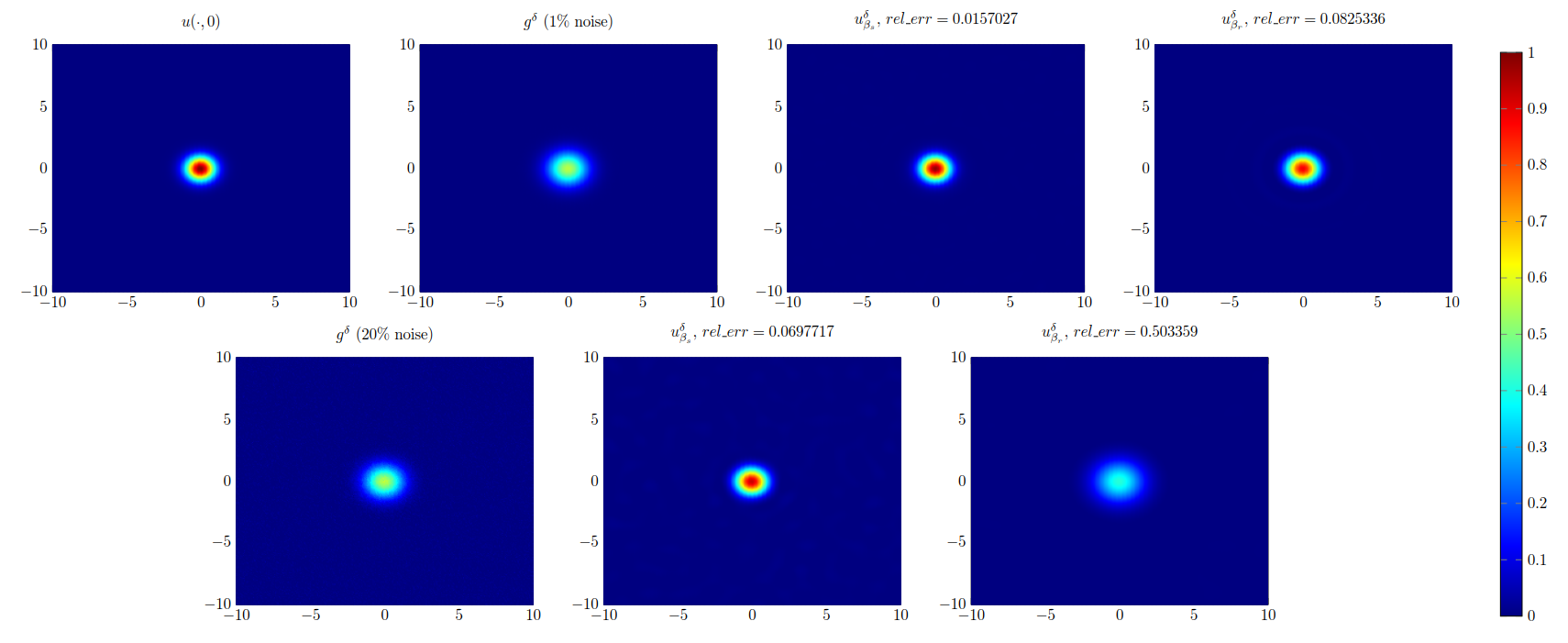} 
\end{center}
\caption{Illustration reconstructed solutions $u_{\beta_s}^\delta$ and $u_{\beta_r}^\delta$ for $1\%$ (first row) and $20\%$ (second row) noise level in Example 1.}
\label{Fig comp sol 1 percent noise level}
\end{figure}

\begin{figure}[h!]
\begin{center}
\includegraphics[scale=0.4]{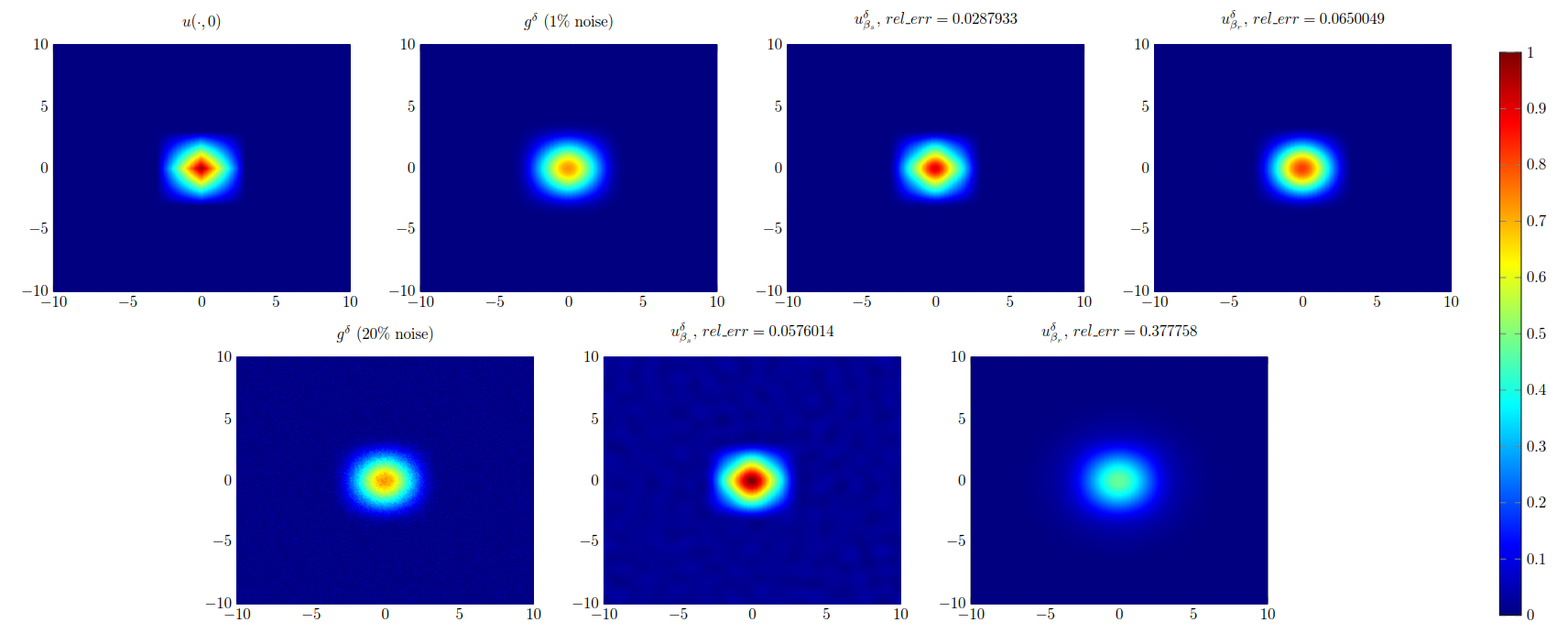} 
\end{center}
\caption{Illustration reconstructed solutions $u_{\beta_s}^\delta$ and $u_{\beta_r}^\delta$ for $1\%$ (first row) and $20\%$ (second row) noise level in Example 2.}
\label{Fig comp sol 2 percent noise level}
\end{figure}

\begin{figure}[h!]
\begin{center}
\includegraphics[scale=0.4]{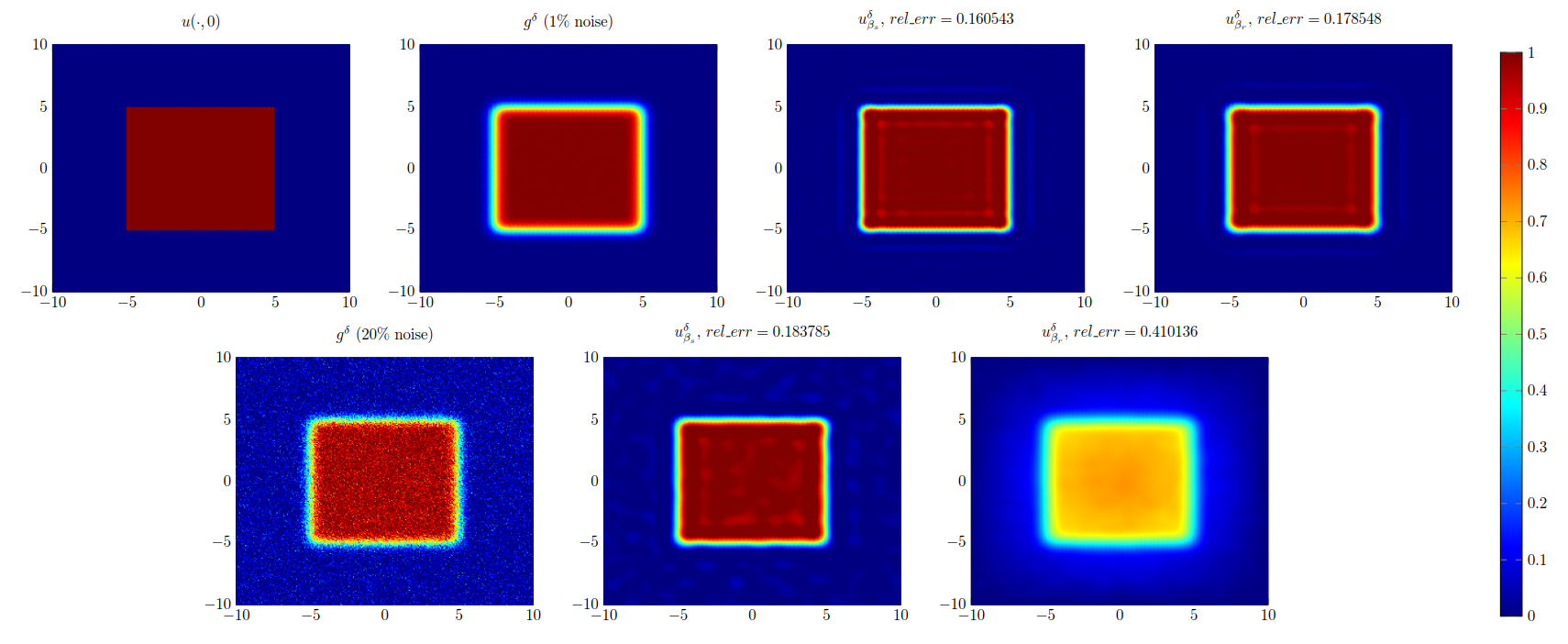} 
\end{center}
\caption{Illustration reconstructed solutions $u_{\beta_s}^\delta$ and $u_{\beta_r}^\delta$ for $1\%$ (first row) and $20\%$ (second row) noise level in Example 3.}
\label{Fig comp sol 3 percent noise level}
\end{figure}

\begin{figure}[h!]
\begin{center}
\includegraphics[scale=0.4]{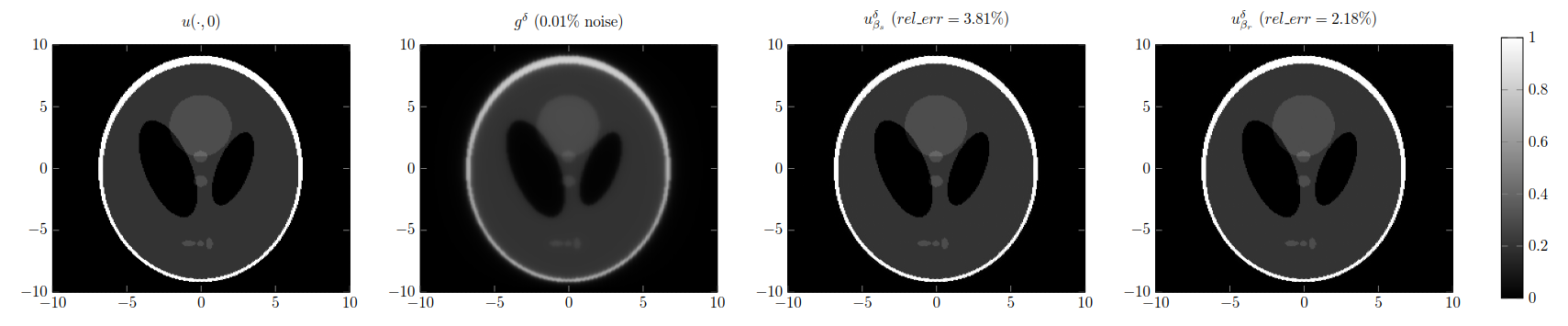} 
\end{center}
\caption{Illustration reconstructed solutions $u_{\beta_s}^\delta$ and $u_{\beta_r}^\delta$ for $0.01\%$ noise level in Example 4.}
\label{Fig comp sol 4 percent noise level}
\end{figure}

Given a reconstructed solution $u_\beta^\delta$, we consider the relative error
$$
rel\_err = \frac{|| u_\beta^\delta- u(\cdot,0)||_2}{||u(\cdot,0)||_2}.
$$
For the selection of the regularization parameter $\beta$, we consider two rules: the a-priori order-optimal rule $\beta_s := c \delta^{1/2s}$ which is independent of the a-priori on the solution $u(\cdot,0)$; and the a-posteriori rule $\beta_r := \beta(\delta,g^\delta,r)$ defined in \eqref{def beta a posteriori rule}. 
For the computation of $\beta_r$, we used Algorithm \ref{Algo beta r} with $\beta_0 =10$ and $q =0.98$.
In Example 1, 2 and 3, we set $c=0.2$. In Example 4, in order to recover irregularity of the initial data $u(\cdot,0)$, we need to applied less regularization, therefore, for example 4, we consider one tenth of the value of $c$ taken for Examples 1,2 and 3. That is, we set $c=0.02$ in Example 4. In all the four examples, for the a-posteriori rule \eqref{def beta a posteriori rule}, we set $r=1$. We recall that with this choice, we obtain sub-optimal rate for Example 1 (see Remark \ref{remark 6}). For example 2 and 3, the order optimality is not guaranteed. This choice of $r=1$ is motivated by simulation observation: in the simulations, this value of $r$ is preferable that smaller values of $r$.

On Figures \ref{Fig comp sol 1 percent noise level} to \ref{Fig comp sol 3 percent noise level}, we illustrate the noisy data $g^\delta$ and the corresponding reconstructed solution $u_{\beta_s}^\delta$ and $u_{\beta_r}^\delta$ for $1\%$ (first row) and $20\%$ (second row) noise level for Example 1 to 3. On these Figures, we can see that even for a $20\%$ noise level, in Examples 1,2 and 3, the reconstructed solution $u_{\beta_s}^\delta$ is still reasonable and exhibit key features of the true solution. Figure \ref{Fig comp sol 4 percent noise level} exhibits reconstructed solutions for the deblurring image example for $0.01\%$ noise level. We recall that in case of image deblurring problem modeled by diffusion process, reasonable reconstruction are obtained only for low level noise (see, e.g. \cite{wang2013total} where blurring is modeled by slow diffusion using time-fractional diffusion equation).

\begin{figure}[h!]
\begin{center}
\includegraphics[scale=0.45]{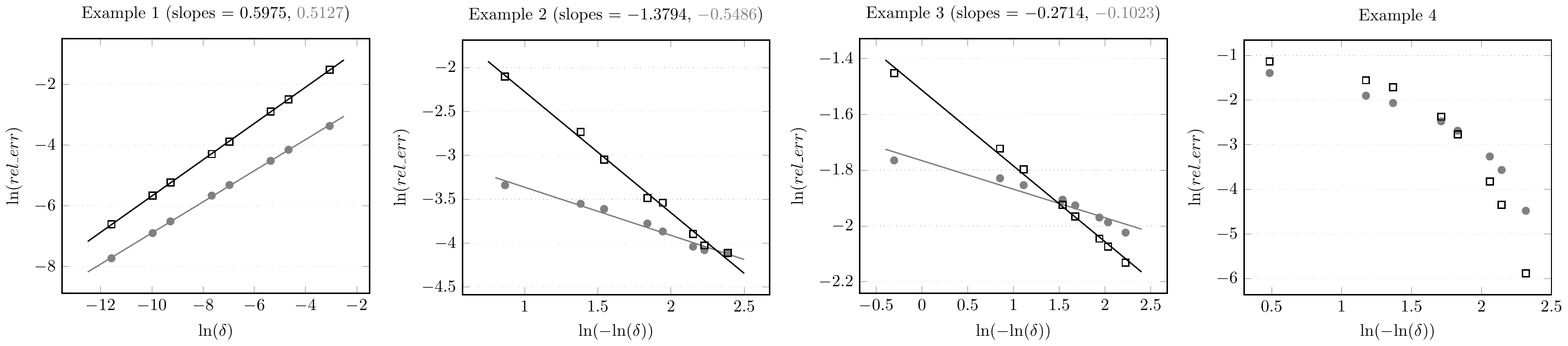} 
\end{center}
\caption{Illustration of numerical rates of convergence for the a-priori rule $\beta_s = c \delta^{1/2s}$ (dark) and a-posteriori rule \eqref{def beta a posteriori rule} with $r=1$ (gray).}
\label{Fig ill num conv rate}
\end{figure}

In order to check the logarithmic convergence rates of the parameter selection rules $\beta_s$ and $\beta_r$, we plot $\ln(rel\_err)$ versus $\ln(-\ln(\delta))$ for various values of $\delta$ for each rule on Figure \ref{Fig ill num conv rate}. We recall that if $rel\_err = \bigO{f_q(\delta^2)}$ as $\delta \to 0$, then the curve $(\ln(-\ln(\delta)),\ln(rel\_err))$ should exhibit a line shape with slope equal to $-q$. Similarly if $rel\_err = \bigO{\delta^q}$ as $\delta \to 0$, then the curve $(\ln(\delta),\ln(rel\_err))$ should exhibit a line shape with slope equal to $q$.

On Figure \ref{Fig ill num conv rate}, in the left plot, we can see that the rate of convergence of error in example 1 is perfectly linear, not logarithmic. This can be explained by the fact that in example, $u(\cdot,0) \in H^p(\Rb^2)$ for all $p>0$. In the second and third plot of Figure \ref{Fig ill num conv rate}, we can see that the reconstructed errors of both rules in example 2 and 3 exhibit logarithmic rates confirming the theoretical logarithmic rate. Moreover, the numerical order $q_{num}$ of the logarithmic rate, though different from the theoretical order $q = p/2$, decreases from Example 2 to Example 3, confirming also the theoretical prediction given the the analytic order $q=p/2$ also decreases from Example 2 to Example 3. In adequation to the theory, the reconstructed error in Example 4 does not exhbit a logarithmic rates, probably due to the fact that in Example 4, $u(\cdot,0)$ does not satisfy smoothness condition \eqref{smoothness cond on u(cdot,0)}. In Example 1,2 and 3, from Figure \eqref{Fig ill num conv rate}, we can see that the a-priori rule $\beta_s$ exhibits higher convergence order compared to the a-posteriori rule \eqref{def beta a posteriori rule}. 

For assessing the numerical stability and convergence of our method, we run a Monte Carlo simulations of $200$ replications of noise term for each example with various noise level. The results are summarized in Tables \ref{Table MC exple 1,2} and \ref{Table MC exple 3,4}. From these Tables , we can confirm the numerical convergence of the method for both rules. This is illustrated by the fact that both the average reconstruction error ($\mathrm{mean}(rel\_err)$) and the average regularization parameters decreases as the noise level decreases. Moreover the very small magnitude of the variance ($\mathrm{var}(rel\_err)$) of the reconstructed error  indicates the numerical stability of the regularization scheme. As predicted according to the regularity, under a same noise level, the quality of reconstruction get better as we move from Example 3 to Example 1.

\begin{table}
\begin{center}
\includegraphics[scale=0.9]{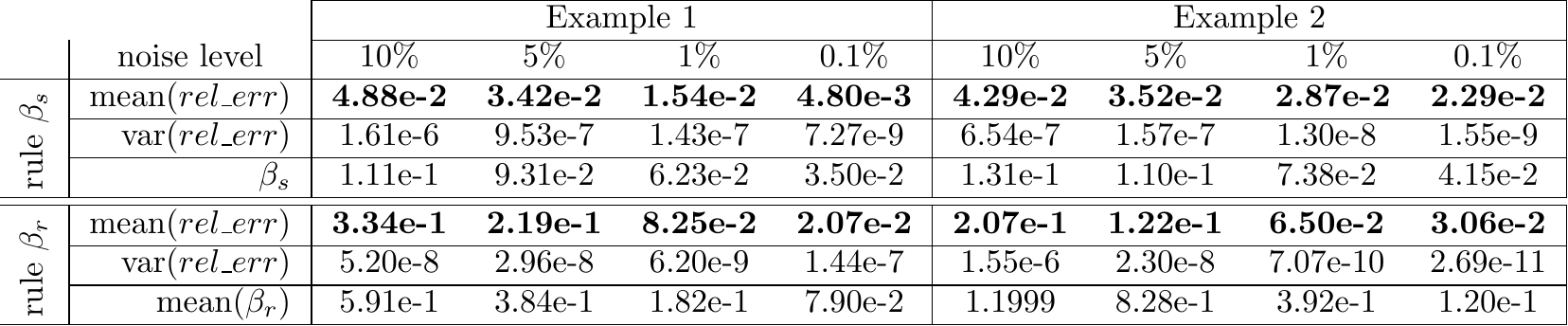}
\end{center}
\caption{Summary Monte Carlo simulation for Example 1 and 2 with 200 sample size.}
\label{Table MC exple 1,2}
\end{table}

\begin{table}
\begin{center}
\includegraphics[scale=0.9]{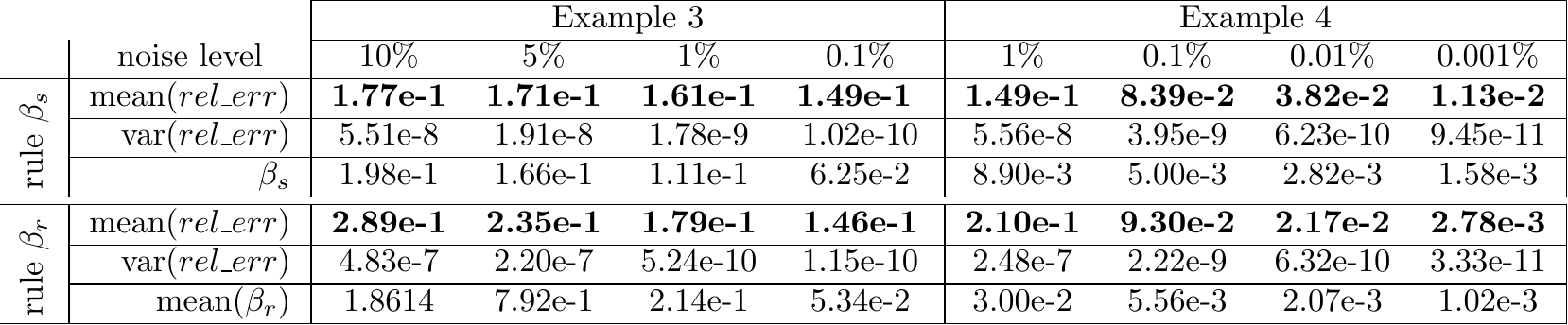} 
\end{center}
\caption{Summary Monte Carlo simulation for Example 3 and 4 with 200 sample size.}
\label{Table MC exple 3,4}
\end{table}

\section{Conclusion}	
We described a simple variational regularization method based on mollification for solving the Backward heat equation with time dependent thermal conductivity and fractional Laplacian in an unbounded domain. We prove order-optimality of the method under Sobolev smoothness condition \eqref{smoothness cond on u(cdot,0)} and derive an order-optimal a-posteriori parameter choice rule based on Morozov principle. We successfully applied the technique for the regularization of three two-dimensional examples together with an example of image deblurring and we were able to obtain reasonable reconstructed solution under high noise level in the data, illustrating so the good potential of the method. Moreover, we numerically confirmed the theoretical convergence rates of the error between the exact solution and the reconstructed solution. We point out that the regularization technique we presented can also be applied to backward heat problems on bounded domains. However, the study of error estimates and convergence rates is much more difficult due to the non-commutation of the mollifier operator $C_\beta$ and the operator $A$ corresponding to the case of bounded domain. 


\section*{Appendix}
The lemma below exhibits some estimates about the logarithmic source function $f_q$ which is repeatedly used in the paper.
\begin{Lemma}
\label{Lemma prop log sour func}
Let $q>0$ and the function $f_q$ defined by $ \forall t \in (0,1]$, $f_q(t) = (-\ln(t))^{-q}$. Then for all $a,b >0$ we have
\begin{equation}
\label{Property log sourc function}
\begin{cases}
\vspace{0.1cm}
\text{if}\,\, \lambda \leq 1, & \forall t \in (0,1), \quad f_q(\lambda t^a) \leq \max \left\lbrace 1,\left(\frac{b}{a}\right)^q \right\rbrace  f_q(\lambda t^b)\\\vspace{0.1cm}
\text{if}\,\, \lambda >1,     & \forall t \in (0, \lambda^{-\frac{2}{a}}), \quad f_q(\lambda t^a) \leq \max \left\lbrace 1, \left( \frac{2b-a}{a}\right)^q \right\rbrace  f_q(\lambda t^b)
\end{cases}
\end{equation}
Moreover, 
\begin{equation}
\label{Property log sourc function 2}
\begin{cases}
\vspace{0.cm}
\text{if}\,\, \lambda \leq 1, & \forall t \in (0,1), \quad f_q(\lambda t) \leq   f_q(t)\\\vspace{0.cm}
\text{if}\,\, \lambda >1,     & \forall t \in (0, \lambda^{-2}), \quad f_q(\lambda t) \leq 2^q   f_q(t).
\end{cases}
\end{equation}
\end{Lemma}
\begin{Proof} We have
\begin{equation}
\label{eqref1}
\forall t \in (0,1), \quad \frac{f_q(\lambda t^a)}{f_q(\lambda t^b)}  = j(x) := \left( \frac{ 1 + a x}{1 + bx} \right)^{-q}, \quad \text{with} \quad x = \frac{\ln t}{\ln \lambda}.
\end{equation}
If $\lambda \in (0,1)$ then $x>0$. For $a \geq b$,  $j(x)$ is obviously less than $1$. For $a < b$ the function $j$ is increasing $j(x)$ is bounded by above on $\Rb^+$ by $\lim_{x \to +\infty} j(x) = (a/b)^{-q}$.\\
Now for $\lambda >1$, and $t \in (0, \lambda^{-\frac{2}{a}})$, we have $x \leq -2/a$. For $a \geq b$, then right hand side in \eqref{eqref1} is bounded above by $1$. For $a < b$, the function $j$ is increasing and thus bounded by above on $(-\infty,-2/a)$ by $j(-2/a) = \left( \frac{a}{2b-a}\right)^{-q}$.
The first inequality in \eqref{Property log sourc function 2} follows readily from the fact the function $f_q$ is increasing on $(0,1)$. For the second inequality, we have
\begin{equation}
\label{eqref2}
\forall t \in (0,1), \quad \frac{f_q(\lambda t)}{f_q(t)}  = k(x) := \left( 1 + y \right)^{-q}, \quad \text{with} \quad y = \frac{\ln \lambda}{\ln t}.
\end{equation}
For $t \leq \lambda^{-2}$, $y>-1/2$, and since the function $k$ in \eqref{eqref2} is decreasing, we deduce that for all $y>-1/2$, $k(y) \leq k(-1/2) = 2^q$ whence the second inequality in \eqref{Property log sourc function 2}.
\end{Proof}

\textbf{Proof of Proposition} \ref{Prop arbitrary slow convergence of C_beta}.
Assume that there exists a function $\nu: \Rb^+ \to \Rb^+$ such that $\lim_{\beta \downarrow 0} \nu(\beta) = 0$ and \eqref{unif bound I-Cbeta} holds. Then for all $f \in L^2(\Rb^n)$ such that $\norm{f}_{L^2} = 1$, using the Parseval identity we get
\begin{equation}
\label{eq 1}
||(1-\hat{\phi}_\beta)\hat{f}||_{L^2} \leq \nu(\beta) \quad i.e. \int_{\Rb^n} \mod{1 - \hat{\phi}(\beta \xi)}^2 \mod{\hat{f}(\xi)}^2 \ud \xi \leq \nu(\beta)^2.
\end{equation}
Since $\hat{\phi}(\xi) \to 0$ as $\mod{\xi} \to \infty$, we get that $\mod{1 - \hat{\phi}(\beta \xi )} \to 1 $ as as $\mod{\xi} \to \infty$. This means that there exists $k_\beta >0$ such that 
\begin{equation}
\label{eq 2}
\mod{\xi} \geq k_\beta \Rightarrow \mod{1 - \hat{\phi}(\beta \xi )} \geq 1/2.
\end{equation}
Now for $k \geq k_\beta +1 $, let 
\begin{equation}
\label{eq 3}
f_k = \Fc^{-1}\left(c \, 1_{B(k e_1,1)} \right),
\end{equation}
where $e_1$ is the first vector in the canonical basis of $\Rb^n$, $B(k e_1,1)$ is the ball of radius $1$ centered at $k\,e_1$ and $c^{2}$ is a normalization coefficient equal to the inverse of the volume of the unit sphere in $\Rb^n$.
We have 
\begin{equation}
\label{eq 4}
\forall \,\xi \in B(k e_1,1), \quad k_\beta \leq k - 1 = \mod{k e_1} - 1 \leq \mod{k e_1} - \mod{\xi - k\,e_1} \leq \mod{\xi} .
\end{equation}
From \eqref{eq 3}, we deduce that 
\begin{equation}
\label{eq 5}
\begin{cases}
\norm{f_k} = 1 \\
\Fc(f_k) = \widehat{f_k} = c\, 1_{B(k e_1,1)}.
\end{cases}
\end{equation}
Applying \eqref{eq 1} with $f_k$ and using \eqref{eq 2},\eqref{eq 4} and \eqref{eq 5} yields
\begin{eqnarray*}
\nu(\beta)^2 \geq \int_{\Rb^n} \mod{1 - \hat{\phi}(\beta \xi)}^2 \mod{\widehat{f_k}(\xi)}^2 \ud \xi  = c^2 \int_{B(k e_1,1)} \mod{1 - \hat{\phi}(\beta \xi )}^2 \ud \xi 
 \geq  \frac{c^2}{4} \int_{B(p e_1,1)} \ud \xi 
 =  \frac{1}{4}
\end{eqnarray*}
By letting $\beta$  goes to $0$, we get that $0 \geq 1/4$. Whence the contradiction. \QED

\textbf{Proof of Lemma} \ref{lemma completion lemma alibaud}. i) follows readily from \eqref{cond on conv kernel phi}.\\
ii) Let $p \in \Rb$ and $f \in H^p(\Rb^n)$,
if $p\leq s$,
\begin{eqnarray*}
 || (I - C_\beta)f||_{L^2}^2 &=& \int_{\Rb^n} |1-\hat{\phi}(\beta \xi)|^{2(1-p/s)}\,  |1-\hat{\phi}(\beta \xi)|^{2p/s} |\hat{f}(\xi)|^2 \ud \xi \quad \text{by Parseval identity}\\
 & \leq & (1 + ||\phi||_{L^1(\Rb^n)})^{2(1-p/s)} \int_{\Rb^n} |1-\hat{\phi}(\beta \xi)|^{2p/s} |\hat{f}(\xi)|^2 \ud \xi  \\
 & = & (1 + ||\phi||_{L^1(\Rb^n)})^{2(1-p/s)} \int_{\Rb^n} |1-\hat{\phi}(\beta \xi/\vert \xi \vert)|^{2p/s} \mod{\frac{|1-\hat{\phi}(\beta \xi)|^2}{|1-\hat{\phi}(\beta \xi/\vert \xi \vert)|^{2}} }^{p/s} |\hat{f}(\xi)|^2 \ud \xi \\
& \leq & (1 + ||\phi||_{L^1(\Rb^n)})^{2(1-p/s)} C_0^{p/s} M_\beta^{p/s}  \int_{\Rb^n} \vert \xi \vert^{2p} |\hat{f}(\xi)|^2 \ud \xi \quad \text{from} \,\, \eqref{def m_beta and M_beta} \,\, \text{and} \,\, \eqref{key estimate alibaud} \\
& \leq & (1 + ||\phi||_{L^1(\Rb^n)})^{2(1-p/s)} C_0^{p/s} M_\beta^{p/s}  \int_{\Rb^n} (1+\vert \xi \vert^{2})^p |\hat{f}(\xi)|^2 \ud \xi \\
& \leq &  C_1 \beta^{2p} ||f||_{H^p}^2 \quad \text{using}\quad i) .
\end{eqnarray*}
For $p>s$,
\begin{eqnarray*}
\label{eqx}
 || (I - C_\beta)f||_{L^2}^2 &=& \int_{\Rb^n} |1-\hat{\phi}(\beta \xi/\vert \xi \vert)|^2 \frac{|1-\hat{\phi}(\beta \xi)|^2}{|1-\hat{\phi}(\beta \xi/\vert \xi \vert)|^2} |\hat{f}(\xi)|^2 \ud \xi  \nonumber \\
& \leq & M_\beta C_0 \int_{\Rb^n} \vert \xi \vert^{2s} |\hat{f}(\xi)|^2 \ud \xi \quad \text{from} \,\, \eqref{def m_beta and M_beta} \,\, \text{and} \,\, \eqref{key estimate alibaud} \\
& \leq  &  C_1 \beta^{2s} ||f||_{H^s}^2 \quad \text{using}\quad i) \nonumber.
\end{eqnarray*}
iii) Let $p \in \Rb$ and $f \in H^{2p}(\Rb^n)$, if $p\leq s$, 
\begin{eqnarray*}
 || (I - C_\beta)^*(I - C_\beta)f||_{L^2}^2 &= &\int_{\Rb^n} |1-\hat{\phi}(\beta \xi)|^{4(1-p/s)}\,  |1-\hat{\phi}(\beta \xi)|^{4p/s} |\hat{f}(\xi)|^2 \ud \xi \quad \text{by Parseval identity} \\
 & \leq &  \tilde{C} \int_{\Rb^n} |1-\hat{\phi}(\beta \xi)|^{4p/s} |\hat{f}(\xi)|^2 \ud \xi \quad \text{with} \quad \tilde{C} = (1 + ||\phi||_{L^1(\Rb^n)})^{4(1-p/s)} \\
 & = & \tilde{C} \int_{\Rb^n} |1-\hat{\phi}(\beta \xi/\vert \xi \vert)|^{4p/s} \mod{\frac{|1-\hat{\phi}(\beta \xi)|^2}{|1-\hat{\phi}(\beta \xi/\vert \xi \vert)|^{2}} }^{2p/s} |\hat{f}(\xi)|^2 \ud \xi \\
& \leq&  \tilde{C} C_0^{2p/s} M_\beta^{2p/s}  \int_{\Rb^n} \vert \xi \vert^{4p} |\hat{f}(\xi)|^2 \ud \xi \quad \text{from} \,\, \eqref{def m_beta and M_beta} \,\, \text{and} \,\, \eqref{key estimate alibaud} \\
& \leq &  C_2 \beta^{4p} ||f||_{H^{2p}(\Rb^n)}^2 \quad \text{using}\quad i) .
\end{eqnarray*}
For $p>s$, 
\begin{eqnarray*}
\label{eqxbis}
 || (I - C_\beta)^*(I - C_\beta)f||_{L^2}^2 &=& \int_{\Rb^n} |1-\hat{\phi}(\beta \xi/\vert \xi \vert)|^4 \mod{\frac{|1-\hat{\phi}(\beta \xi)|^2}{|1-\hat{\phi}(\beta \xi/\vert \xi \vert)|^2}}^2  |\hat{f}(\xi)|^2 \ud \xi  \nonumber \\
& \leq & M_\beta^2 C_0^2 \int_{\Rb^n} \vert \xi \vert^{4s} |\hat{f}(\xi)|^2 \ud \xi \quad \text{from} \,\, \eqref{def m_beta and M_beta} \,\, \text{and} \,\, \eqref{key estimate alibaud} \\
& \leq  &  C_2 \beta^{4s} ||f||_{H^{2s}}^2 \quad \text{using}\quad i) \nonumber. \qquad \qquad \qquad \qquad \qquad \qquad \qquad \qquad \QED
\end{eqnarray*}

\bibliographystyle{abbrv}
\bibliography{References_paper_TAO}
\end{document}